\documentclass[12pt, a4paper]{amsart}
\input cyracc.def
\newfam\cyrfam
\font\tencyr=wncyr10
\font\sevencyr=wncyr7

\def\cyr{\fam\cyrfam\tencyr\cyracc}
\textfont\cyrfam=\tencyr \scriptfont\cyrfam=\sevencyr
 \scriptscriptfont\cyrfam=\sevencyr

\font\tencyr=wncyr10

\newcommand {\Cee} {{\mathbb C}}
\newcommand {\Kee} {{\mathbb K}}
\newcommand {\Qee} {{\mathbb Q}}
\newcommand {\Ree} {{\mathbb R}}
\newcommand {\Zee} {{\mathbb Z}}

\newcommand {\cal} {\mathcal}
\newcommand {\cD} {{\cal D}}

\newcommand {\fag} {{\mathfrak{ag}}}

\newcommand {\fc} {{\mathfrak{c}}}
\newcommand {\fp} {{\mathfrak{p}}}
\newcommand {\fs} {{\mathfrak{s}}}
\newcommand {\fz} {{\mathfrak{z}}}

\newcommand {\fder} {{\mathfrak{der}}} %

\newcommand {\fe} {{\mathfrak{e}}}

\newcommand {\ff} {{\mathfrak{f}}}

\newcommand {\fg} {{\mathfrak{g}}} %
\newcommand {\fgl} {{\mathfrak{gl}}} %
\newcommand {\fh} {{\mathfrak{h}}}
\newcommand {\fhei} {{\mathfrak{hei}}}
\newcommand {\fk} {{\mathfrak{k}}}
\newcommand {\fl} {{\mathfrak{l}}}

\newcommand {\fm} {{\mathfrak{m}}}
\newcommand {\fn} {{\mathfrak{n}}}

\newcommand {\fo} {{\mathfrak{o}}}
\newcommand {\fpo} {{\mathfrak{po}}}
\newcommand {\fsl} {{\mathfrak{sl}}}

\newcommand {\fsp} {{\mathfrak{sp}}}

\newcommand {\fsvect} {{\mathfrak{svect}}}
\newcommand {\fvect} {{\mathfrak{vect}}} %
\newcommand {\fv} {{\mathfrak{v}}}

\def \opname#1#2%
 {\expandafter\newcommand \csname #1\endcsname {{\mathop{#2}\nolimits}}}

\newcommand{\rmname}[1]
 {\expandafter\newcommand \csname #1\endcsname {{\operatorname{#1}}}}

\newcommand{\rmnameii}[2]
 {\expandafter\newcommand \csname #1\endcsname {{\operatorname{#2}}}}

\rmname{card} \rmname{gr} \rmname{Lie} \rmname{rk} \rmname{Hom}
\rmname{Span}

\rmnameii {IM} {Im}
\rmnameii {RE} {Re}

\opname{GL} {{G\hspace{-0.3ex}L}}
\opname{SL} {{S\hspace{-0.3ex}L}}

\newcommand {\tto} {\longrightarrow}
\newcommand {\pder}[1] {{\frac{\partial}{\partial {#1}}}}
\newcommand {\pderf}[2] {{\frac{\partial {#1}}{\partial {#2}}}}

\newcommand {\bcdot} {\mathbin{\hbox{\raise.4ex\hbox{\bf.}}}}

\def \nopoint#1#2{}

\newcommand {\noqed}
 {\renewcommand {\qed} {}}

\newcommand {\secno} {}
\newcommand {\ssecfont} {\normalfont\bf}

\newtheorem*{Theorem}{\secno Theorem}

\newtheorem*{Corollary}{\secno Corollary}
\newtheorem*{Statement}{\secno Statement}

\newenvironment {th*}[1]
 {\gdef\thname{#1} \begin{thn}}%
 {\end{thn}}
\newtheorem*{thn} {\thname}

\theoremstyle{definition}

\newenvironment {ex*}[1]
 {\gdef\thname{#1} \begin{exn}}%
 {\end{exn}}
\newtheorem*{exn}{\thname}

\theoremstyle{remark}

\newtheorem*{Remark}{\secno Remark}

\newenvironment {rem*}[1]
 {\gdef\thname{#1} \begin{remn}}%
 {\end{remn}}
\newtheorem*{remn}{\thname}

\newcommand {\ssec}{\subsection*}

\newcommand {\ssbegin}[2]
 {\def \secno {\gdef \secno {}{\ssecfont #1. }}%
 \begin{#2}}
\begin{document}

\title{The nonholonomic Riemann and Weyl tensors for flag manifolds}

\author{Pavel Grozman${}^1$, Dimitry Leites${}^2$}

\address{${}^1$Equa Simulation AB, Stockholm, Sweden; pavel@rixetele.com\\
${}^2$MPIMiS, Inselstr. 22, DE-04103 Leipzig, Germany \\
on leave from Department of Mathematics, University of Stockholm,
Roslagsv. 101, Kr\"aftriket hus 6, SE-104 05 Stockholm, Sweden\\
mleites@math.su.se}



\begin{abstract} On any manifold, any non-degenerate symmetric 2-form (metric)
and any skew-symmetric (differential) form $\omega$ can be reduced to a
canonical form at any point, but not in any neighborhood: the respective
obstructions being the Riemannian tensor and $d\omega$. The obstructions to
flatness (to reducibility to a canonical form) are well-known
for any $G$-structure, not only for Riemannian or symplectic structures.

For the manifold with a nonholonomic structure
(nonintegrable distribution), the general notions of flatness and
obstructions to it, though of huge interest (e.g., in supergravity) were not
known until recently, though particular cases were known for more than a century
(e.g., any contact structure is ``flat'': it can always be reduced, locally, to
a canonical form).

We give a general definition of the {\it nonholonomic} analogs of
the Riemann and Weyl tensors. With the help of Premet's theorems
and a package {\bf SuperLie} we calculate these tensors for the
particular case of flag varieties associated with each maximal
(and several other) parabolic subalgebra of each simple Lie
algebra. We also compute obstructions to flatness of the
$G(2)$-structure and its nonholonomic super counterpart.
\end{abstract}

\thanks{DL thankfully acknowledges help of A.~Premet, I.~Shchepochkina and
J.~Bernstein and financial support of IHES, F\"ors\"akringskassan
and MPIMiS}

\maketitle

\begin{quote}\rightline{To Antony Joseph}\end{quote}

\section*{Introduction}

H.~Hertz \cite{H} coined the term {\it nonholonomic} during his
attempts to geometrically describe motions in such a way as to
exorcize the concept of ``force''. A manifold (phase space) is
said to be {\it nonholonomic} if endowed with a nonintegrable
distribution (a subbundle of the tangent bundle). A simplest
example of a nonholonomic dynamical system is given by a body
rolling without gliding over another body. Among various images
that spring to mind, the simplest is a ball on a rough plane
(\cite{Poi}) or a bike on asphalt. At the tangency point of the
wheel with asphalt, the velocity of the wheel is zero. (This is a
linear constraint. We will not consider here more general
non-linear constraints.) A famous theorem of Frobenius gives
criteria of local integrability of the distribution: its sections
should form a Lie algebra.

For a historical review of nonholonomic systems, see \cite{VG2}
and a very interesting paper by Vershik \cite{V} with first
rigorous mathematical formulations of nonholonomic geometry and
indications to applications to various, partly unexpected at that
time, areas (like optimal control or macro-economics, where {\it
nonlinear} constraints are also natural, cf. \cite{AS}, \cite{Bl},
\cite{S}); recent book by Kozlov \cite{Koz} is extremely
instructive. In \cite{V}, Vershik summarizes about 100 years of
studies of nonholonomic geometry (Hertz, Carath\'eodory, Vr\u
anceanu, Wagner, Schouten, Faddeev, Griffiths, Godbillon; now
MathSciNet returns thousands entries for ``nonholonomic'' and its
synonyms (anholonomic, ``sub-Riemannian'', ``autoparallel'') and
particular cases leading to nonholonomic constraints (``Finsler'',
``cat's problem''). There seems to be ``more", actually,
nonholonomic dynamical systems than holonomic ones.

A relatively new theory, ``supergravity'' (the theory embodying Einstein's dream
of a Unified Field Theory), also deals exclusively with nonholonomic structures,
albeit on {\it super}manifolds.)

At the end of \cite{V} Vershik summarized futile
attempts of the researchers to define an analog of the Riemann tensor for the
general nonholonomic manifold in a conjecture that ``though known in some cases,
it is probably impossible to define such a general analog''.

However, in 1989, during his stay at IAS, DL gave such a general definition
and lectured on it at various schools and conferences (ICTP, Euler Math.
Inst., JINR, etc.), see \cite{L}, \cite{LP}; later we applied it to supergravity
\cite{GL1}.

{\bf Main results of this paper are (1) elucidation of this general definition
of the nonholonomic counterparts of the Riemann tensor and its conformal, Weyl,
analog; (2) Premet's theorems that facilitate computation of these tensors in
some cases (for flag varieties $G/P$, where $G$ is a simple Lie group and $P$
is its parabolic subalgebra); (3) computation of these tensors in some of these cases.}

Tanaka \cite{T} tackled the same, actually, problem for totally
different reasons (in \cite{T} and even in more recent \cite{Y}
and \cite{YY} even the word \lq\lq nonholonomic" is never used).
Tanaka's results (especially their lucid exposition in Yamaguchi's
paper \cite{Y}) are easier to understand than the first attempts
(by Schouten, Wagner, see \cite{DG}) because, after some
experiments, he used the hieroglyphics of Lie algebra cohomology
which are much more graphic than coordinate tensor notations.
Tanaka's tensor coincides with the one we suggest.

We illustrate the main definitions by computing the nonholonomic analogs of the
Riemann and Weyl tensor in several particular cases --- the simplest analogs of
\lq\lq classical domains". In doing so we rely on Premet's theorems and a {\it
Mathematica}-based package {\bf SuperLie}.

Computations of nonholonomic analogs of Riemann tensor are rather
difficult technically and the rare examples of works with actually
computed {\bf results} are \cite{C1}--\cite{C5}, \cite{GIOS},
\cite{HH}, \cite{Y, YY, EKMR, Ta} and refs. therein. In non-super
setting, they used Tanaka's definition of nonholonomic Riemann
tensor, identical to ours, but lack Premet's theorems and {\bf
SuperLie} and so could not compute as much as anybody is able to
now with their help.\footnote{In 2000, S.~Vacaru informed us of
his and mathematician's from Vr\u anceanu's school definitions
partly summarized in \cite{Va6} and refs. therein. It is not easy
to see through the forest of non-invariant cumbrous tensor
expressions that a number of components is lacking in \cite{Va6}
as well as in \cite{DG}, as compared with Tanaka's or our
definitions.}

\ssec{0.1. General description of classical tensors and our examples}
In mid-1970s, Gindikin formulated a problem of local characterization
of compact Hermitian symmetric domains $X=S/P$, where $S$ is a simple Lie group and
$P$ its parabolic subgroup. Goncharov solved this
problem \cite{Go} having considered the fields of certain quadratic
cones and having computed the {\bf structure functions} (obstructions to
flatness) of the corresponding $G$-structures, where $G$ is the
Levi (reductive) part of $P$.

\ssbegin{0.1.1}{Examples} Let the ground field be $\Cee$.

1) For $S=O(n+2)$ and $G=CO(n)=O(n)\times \Cee^\times$, the structure
functions were
known; they constitute the Weyl tensor --- the conformally invariant
part of the Riemann curvature tensor.

2) For $S=SL(n+m)$ and $G=S(GL(n)\times GL(m))$, the structure functions
are obstructions to integrability of multidimensional analogs of
Penrose's $\alpha$- and $\beta$-planes on the Grassmannian
$Gr_n^{n+m}$ (Penrose considered
$Gr_2^{4}$). \end{Examples}

Not any simple complex Lie group $S$ and its subgroup $P$ can form
a classical domain: $S$ is any but $G(2)$\footnote{We denote the
exceptional groups and their Lie algebras in the same way as the
serial ones, like $SL(n)$; we thus avoid confusing $\fg(2)$ with
the second component $\fg_2$ of a $\Zee$-grading of a Lie algebra
$\fg$.}, $F(4)$ and $E(8)$ and $P=P_{i}$ is a {\it maximal}
parabolic subgroup generated by all Chevalley generators of $S$,
but {\bf one} ($i$th), say, negative. The group $P$, or which is
the same, the $i$th Chevalley generator of $S$ (in what follows
referred to as {\it selected}) can not be arbitrary, either. To
describe the admissible $P$'s, let us label the nodes of the
Dynkin graph of $S$ with the coefficients of the maximal root
expressed in terms of simple roots. The selected generator may
only correspond to the vertex with label 1 on the Dynkin graph.

It is natural to consider the following problem:
$$
\begin{array}{c}
 \text{{\bf For any simple Lie group $S$,}}\\
\text{{\bf fix an arbitrary $\Zee$-grading of its Lie algebra $\fs=\text{Lie}(S)$.}}\\
 \text{{\bf For any subgroup
 $P\subset S$, generated by nonnegative elements of
$\fs$, }}\\
 \text{{\bf what are the analogs of the
Goncharov conformal structure}}\\
\text{{\bf and the corresponding analogs of Riemann and projective
structures,}}\\
\text{{\bf which of these structures are flat, }}\\
\text{{\bf and what are the obstructions to their flatness?}}
\end{array}\eqno{(0.1)}
$$

\ssbegin{0.1.2}{Remark} The adjective \lq\lq arbitrary"
($\Zee$-grading of $\fs$) in (0.1) appeared thanks to J.~Bernstein
who reminded us that parabolic subgroups are a particular case of
such gradings. All $\Zee$-gradings are obtained by setting $\deg
X_i^{\pm}= \pm k_i$, where $k_i\in\Zee$, for the Chevalley
generators $X_i^{\pm}$ and parabolic subgroups appear if $k_i\geq
0$ for all $i$. Recently Kostant \cite{K} considered an analog of
the Borel-Weil-Bott (BWB) theorem --- one of our main tools ---
for the non-parabolic case, but the answer is not yet as algebraic
as we need, so having answered questions (0.1) in full generality
we calculate the nonholonomic invariants for parabolic subgroups
only. \end{Remark}

Modern descriptions of structure functions is usually given in
terms of the {\it Spencer cohomology}, cf. \cite{St} (we will
recall all definitions needed in (0.1) and (0.2) in due course).
Goncharov expressed the structure functions as tensors taking
values in the vector bundle over $X=G/P$, whose fibers at every
point $x\in X$ are isomorphic to each other and to
$$
\renewcommand{\arraystretch}{1.3}
\begin{array}{l}
 H^2(\fg_{-1}; (\fg_{-1}, \fg_{0})_*), \;\text{ where
$\fg_{0}=$ Lie$(G)$, \;\; $\fg_{-1}=T_xX$, }\\
\text{and where
$(\fg_{-1}, \fg_{0})_*=\mathop{\oplus}\limits_{i\geq -1}\fg_i$ is the Cartan
prolong of $(\fg_{-1}, \fg_{0})$.}
\end{array}\eqno{(0.2)}
$$
The conventional representation of the structure functions as bigraded Spencer
cohomology $H^{k, 2}$ can be recovered any time as the homogeneous
degree $k$ component of $H^2(\fg_{-1}; (\fg_{-1}, \fg_{0})_*)$
corresponding to the $\Zee$-grading of $(\fg_{-1}, \fg_{0})_*$.

At about the same time Goncharov got his result, physicists trying
to write down various supergravity equations (for standard or
``exotic'' $N$-SUGRAs, see \cite{WB}, \cite{MaG}, \cite{GIOS}, \cite{HH}) bumped
into the same
problem $(0.1)$ with the supergroup $S=\SL(4|N)$ for $N\leq 8$ and $P$
generated by all the (analogs of the) Chevalley generators of $G$ but
{\bf two}. The corresponding coset superspace $X$ is a flag
supervariety and the difficulties with SUGRAs spoken about, e.g., in
\cite{WB} (\lq\lq we do not know how to define the analog of the Riemann
tensor for $N>2$'', in other words: we do not know what might stand in the
left-hand sides of the SUGRA equations), were caused not by a {\it super}
nature of $X$ but by its {\it nonholonomic} nature.

Shchepochkina introduced nonholonomic generalizations $(\fg_{-}, \fg_{0})_*$ of
Cartan prolongation $(\fg_{-1}, \fg_{0})_*$ for needs of our classification of
simple infinite dimensional Lie superalgebras
of vector fields (\cite{LSh}). She introduced them (together with
several new types of prolongation, e.g., {\it partial} prolongation) in
\cite{Sh1}, \cite{Sh2}, \cite{Sh14}. These generalizations are precisely what is
needed to define the nonholonomic analog of the Weyl and Riemann tensors in the
general case.

Observe that our nonholonomic invariants, though natural analogs of
the curvature and torsion tensors, do not coincide on nonholonomic
manifolds with the classical ones and bearing the same name. Indeed,
on any nonholonomic manifold, there is, by definition, a nonzero
classical torsion (the Frobenius form that to a pair of sections of the
distribution assigns their bracket) while, for example, every contact manifold
is flat in {\it our}
sense. To avoid confusion, we should always add adjective
``nonholonomic'' for the invariants introduced below. Since this is
too long, we will briefly say {\it N-curvature} tensor and specify its degree
(=the order of the structure function) if needed; to require vanishing of the
torsion is analogous of imposing Wess-Zumino constraints \cite{WB}.

The main thing is to answer the questions (0.1). Having done this
(having given appropriate definitions in the general case of
manifolds with nonholonomic structure) we explicitly compute the
analogs of (0.2) --- the space of nonholonomic structure functions
--- possible values of the nonholonomic versions of the Weyl and
Riemann tensors. We do so for the simplest nonholonomic flag
manifolds of the form $S/P$ with {\bf one} selected Chevalley
generator. In most of our cases $(\fg_{-}, \fg_{0})_*=\fs$, the
Lie algebra of $S$, and therefore we can apply the Borel-Weil-Bott
(BWB) theorem (reproduced below; for a nice review, see
\cite{Wo}). If $(\fg_{-}, \fg_{0})_*$ strictly contains $\fs$, we
consider the values of cocycles in $\fs$ as well as in $(\fg_{-},
\fg_{0})_*$.

We cite Premet's theorems that show how to compute the N-Weyl and N-Riemann
tensors
and use the theorems to get an explicit answer.

The implicit form of the answer in \cite{Go} hides phenomena
manifest if the answer is explicit, as in \cite{LPS}, where,
thanks to an explicit form of the answer we suggested some analogs
of Einstein equations (EE) for certain Grassmannians. For the
cases we consider here, a phenomenon similar to that observed in
\cite{LPS} is manifest, e.g., for the nodes at the base of the
forks in $\fe(6)$ and $\fo(8)$. We intend to consider the related
analogs of EE elsewhere.

We illustrate usefulness of computer-aided study by using {\bf
SuperLie} to compute the structure functions for the
$G(2)$-structure, so popular lately, cf. \cite{AW, B, FG}. {\bf
SuperLie} already proved useful in many instances (see \cite{GL}),
and is indispensable for Lie superalgebras: for practically all of
them, there exists nothing as neat as the BWB theorem (\cite{PS}).
We also apply {\bf SuperLie} to compute the structure functions
for a super version of the $G(2)$-structure on the projective
superspace $\Cee P^{1, 7}$ with a nonholonomic distribution.

\section*{\S 1. Structure functions of $G$-structures}

Let $M^n$ be a manifold over a field $\Kee$. Let $FM$ be the frame
bundle over $M$, i.e., the principal $GL(n)$-bundle. Let $G\subset
GL(n)$ be a Lie group. A $G$-{\it structure} on $M$ is a reduction
of the principal $GL(n)$-bundle to the principal $G$-bundle.
Another formulation is more understandable: a $G$-{\it structure}
is a selection of transition functions from one coordinate patch
to another so that they belong to $G$ for every intersecting pair
of patches.

Thus, in the definition of $G$-structure the following characters participate:
$M^n$ and two vector bundles over it: $TM$ and $FM$ and the two groups $G\subset
GL(n)$ both acting in each fiber of each bundle.

The simplest $G$-structure is the {\it flat} $G$-structure defined as
follows. For a model manifold with the flat $G$-structure we take $V=\Cee^n$
with a fixed frame. The key moment is identification of the tangent spaces
$T_vV$ at distinct points $v$. This is performed by means of parallel
translations along $v$. This means that we consider $V$ as a commutative Lie
group and identify the tangent spaces to it at various points with its Lie
algebra, $\fv$. Thanks to commutativity:
$$
\begin{array}{l}
\fv\text{ can be naturally identified with $V$ itself;}\\
\text{it does not matter whether we use left or right translations. }
\end{array}\eqno{(1.1)}
$$

In this way, we get a fixed frame in every $T_vV$. The {\it flat
$G$-structure} is the bundle over $V$ whose fiber over $v\in V$
consists of all frames obtained from the fixed one under the
$G$-action. In textbooks on differential geometry (e.g., in
\cite{St}), the obstructions to identification of the $k$th
infinitesimal neighborhood of a point $m\in M$ on a manifold $M$
with $G$-structure with the $k$th infinitesimal neighborhood of a
point of the manifold $V$ with the above flat $G$-structure are
called {\it structure functions} of order $k$.

To precisely describe the structure functions, set
$$
\fg_{-1} = T_mM,\; \; \fg_0 = \fg = \Lie (G).
$$
Recall that, for any (finite dimensional) vector space $V$, we
have
$$
\Hom(V, \Hom(V,\ldots, \Hom(V,V)\ldots))
\simeq L^{i}(V, V, \ldots, V; V),
$$
where $L^{i}$ is the space of $i$-linear maps and we have
$(i+1)$-many $V$'s on both sides. Now, we recursively define, for
any $i > 0$:
$$
\renewcommand{\arraystretch}{1.4}
\begin{array}{l}
 \fg_i = \{X\in \Hom(\fg_{-1}, \fg_{i-1})\mid X(v_1)(v_2, v_3, ..., v_{i+1})
 =
X(v_2)(v_1, v_3, ..., v_{i+1})\\
\text{ where $v_1, \dots, v_{i+1}\in \fg_{-1}\}$}. \end{array}
$$
Let the $\fg_0$-module $\fg_{-1}$ be faithful. Then, clearly,
$$
(\fg_{-1}, \fg_{0})_*\subset \fvect (n) = \fder~ \Cee[x_1, \ldots
, x_n ]], \; \text{ where}\; n = \dim~ \fg_{-1}.
$$
It is subject to an easy verification that the Lie algebra
structure on $\fvect (n)$ induces same on $(\fg_{-1}, \fg_{0})_*$.
(It is also easy to see that  even if $\fg_{-1}$ is not a faithful
$\fg_0$-module $(\fg_{-1}, \fg_{0})_*$ is a Lie algebra, but can
not be embedded into $\fvect(\fg_{-1}^*)$.) The Lie algebra
$(\fg_{-1}, \fg_{0})_*$ will be called the {\it Cartan's prolong}
(the result of Cartan's {\it prolongation}) of the pair
$(\fg_{-1}, \fg_{0})$. The Cartan prolong is the Lie algebra of
symmetries of the $G$-structure in the space $T_{m}M$.

Let $E^i$ be the operator of the $i$th exterior power, $V^*$ the
dual of $V$. Set
$$
C^{k, s}_{(\fg_{-1}, \fg_0)_*} =
\fg_{k-s}\otimes E^s(\fg_{-1}^*).
$$
 The differential $\partial
_s: C^{k, s}_{(\fg_{-1}, \fg_0)_*}\tto C^{k, s+1}_{(\fg_{-1},
\fg_0)_*}$ is given by (as usual, the slot with the hatted variable is
to be ignored):
$$
(\partial_sf)(v_1, \ldots , v_{s+1}) = \sum _i (-1)^i[f(v_1, \ldots ,
\widehat
{v_{s+1-i}}, \ldots , v_{s+1}), v_{s+1-i}]
$$
for any $v_1, \ldots , v_{s+1}\in \fg_{-1}$. As expected, $\partial
_s\partial _{s+1} = 0$. The homology of this bicomplex is called {\it
Spencer cohomology} of the pair $(\fg_{-1}, \fg_0)$ and denoted by
$H^{k, s}_{(\fg_{-1}, \fg_0)_*}$.

\ssbegin{1.1}{Proposition}[\cite{St}] The order $k$ structure
functions of the $G$-structure --- obstructions to identification
of the $k$th infinitesimal neighborhood of the point in a manifold
with a flat $G$-structure with that at a given point $m\in M$ ---
span, for every $m$, the space $H^{k, 2}_{(\fg_{-1}, \fg_0)_*}$.
These obstructions are defined provided obstructions of lesser
orders vanish.
\end{Proposition}

\ssbegin{1.2}{Example} All structure functions of any
$GL(n)$-structure vanish identically, so all $GL(n)$-structures
are locally equivalent, in particular, locally flat. Indeed: by a
theorem of Serre (\cite{St}) $H^2(V; (V, \fgl(V))_*)=0$ .
\end{Example}

Clearly, the order of the structure functions of a given $G$-structure
may run 1 to $N+2$ (or 1 to $\infty$ if $N=\infty$), where $(\fg_{-1},
\fg_0)_*=\mathop{\oplus}\limits_{i=-1}^N \fg_i$.

\ssbegin{1.3}{Example} Let $\fg_0=\fc\fo(V):=\fo(V)\oplus\Cee z$ be the Lie
algebra of conformal transformations, $\fg_{-1}=V$, $\dim V=n$. For $n=2$, let
$V=V_1\oplus V_2$ with basis $\partial_x$ and $\partial_y$ and let $\fo(V):=\Cee
(x\partial_x-y\partial_y)$. Then ({\it Liouville's theorem}, \cite{St})
\footnotesize
$$
(V, \fc\fo(V))_*= \begin{cases}
\fvect(V^*)&\text{for $n=1$},\\
\fvect(V_1^*)\oplus\fvect(V_2^*)&\text{for $n=2$},\\
V\oplus\fc\fo(V)\oplus V^*\simeq \fo(n+2)&\text{for $n>2$},\end{cases}
\quad
(V, \fo(V))_*= \begin{cases}V&\text{for $n=1$},\\
V\oplus \fo(V)&\text{for $n\geq 2$ }.
\end{cases}
$$\normalsize
The values of the Riemann tensor on any $n$-dimensional Riemannian
manifold belong to $H^{2, 2}_{(V, \fo(V))_*}$
whereas $H^{1, 2}_{(V, \fo(V))_*}=0$.

The fact that $H^{1,
2}_{(V, \fo(V))_*}=0$ (no torsion) is usually referred to as (a part of) the {\it Levi-Civita
theorem}. It implies that, in the Taylor series expansion of the metric at some
point (here $\eta$ is the canonical form; $x$ is the vector of coordinates, so
$x^2$ is the vector of pairs of coordinates, etc.),
$$
g(x)=\eta+ s_1x+ s_2x^2+ s_3x^3+\ldots
$$
the term $s_1$ can be eliminated by a choice of coordinates. Since there are no
structure functions of orders $>2$, all the $s_i$ with $i\geq 2$ only depend on the
Riemann tensor.
\end{Example}

\ssbegin{1.4}{Remark} (cf. \cite{Go}.) Let $H^s_{k}$ be the degree $k$ component
of $H^s(\fg_{-1}; (\fg_{-1},
\fg_0)_*)$ with respect to the $\Zee$-grading induced by the $\Zee$-grading
of $(\fg_{-1},
\fg_0)_*$. Clearly,
$H^{k, s}_{(\fg_{-1}, \fg_0)_*}=H^s_{k}$, so
$$
\mathop{\oplus}\limits_kH^{k, s}_{(\fg_{-1}, \fg_0)_*}= H^s(\fg_{-1};
(\fg_{-1}, \fg_0)_*).
$$
\end{Remark}

This remark considerably simplifies calculations, in particular, if
the Lie algebra $(\fg_{-1}, \fg_0)_*$ is simple and finite
dimensional, we can apply the BWB theorem. In the nonholonomic case
considered in what follows we apply the remark to give a compact
definition\footnote{Cf. with the problems encountered in the pioneer
papers \cite{T}, where only $d=2$ is considered. Wagner's tensors (for any
$d$) look even more horrible, see \cite{DG}, \cite{Va6}.} of structure
functions. We can recover the bigrading at any moment but to
work with just one grading is much simpler.

\section*{\S 2. Structure functions of nonholonomic structures}

To embrace contact-like structures, we have to slightly generalize
the notion of Cartan prolongation: with the tangent bundle over every
nonholonomic manifold there is naturally associated a bundle of graded
nilpotent Lie algebras, cf. \cite{VG}, \cite{M}. For example, for any odd
dimensional manifolds with a contact structure, this is a bundle of
Heisenberg Lie algebras.

\ssec{2.1. Nonholonomic manifolds (\cite{VG, VG2}). Nonholonomic
manifolds. Tanaka-Shchepochkina prolongs} Let $M^n$ be an
$n$-di\-men\-si\-o\-nal manifold with a nonintegrable distribution
$\cD$. Let
$$
\cD= \cD_{-1}\subset \cD_{-2} \subset \cD_{-3} \dots \subset
\cD_{-d}
$$
be the sequence of strict inclusions, where the fiber of
$\cD_{-i}$ at a point $x\in M$ is
$$
\cD_{-i+1}(x)+ [\cD_{-1}, \cD_{-i+1}](x)
$$
(here $[\cD_{-1}, \cD_{-i-1}]=\Span\left([X, Y]\mid X\in
\Gamma(\cD_{-1}), Y\in \Gamma(\cD_{-i-1})\right)$) and $d$ is the
least number such that
$$
\cD_{-d}(x)+[\cD_{-1}, \cD_{-d}](x) = \cD_{-d}(x).
$$
In case $\cD_{-d} = TM$ the distribution is called {\it completely
nonholonomic}. The number $d = d(M)$ is called the {\it
nonholonomicity degree}. A manifold $M$ with a distribution $\cD$
on it will be referred to as {\it nonholonomic} one if $d(M)\neq
1$. Let
$$
n_i(x) = \dim \cD_{-i}(x);
\qquad n_0(x)=0; \qquad n_d(x)=n-n_{d-1}.\eqno{(2.1)}
$$

The distribution $\cD$ is said to be {\it regular} if all the
dimensions $n_i$ are constants on $M$. We will only consider
regular, completely nonholonomic distributions, and, moreover,
satisfying certain transitivity condition (\ref{tr}) introduced
below.

To the tangent bundle over a nonholonomic manifold $(M, \cD)$ we
assign a bundle of $\Zee$-graded nilpotent Lie algebras as
follows. Fix a point $pt\in M$. The usual adic filtration by
powers of the maximal ideal $\fm:=\fm_{pt}$ consisting of
functions that vanish at $pt$ should be modified because distinct
coordinates may have distinct \lq\lq degrees". The distribution
$\cD$ induces the following filtration in $\fm$:
$$
\renewcommand{\arraystretch}{1.4}
\begin{array}{ll}
\fm_k=&\{f\in\fm\mid X_1^{a_1}\ldots X_n^{a_n}(f)=0\;\text{ for
any $X_1, \dots,X_{n_1}\in \Gamma(\cD_{-1})$,  }\\
&  \text{$X_{n_1+1}, \dots, X_{n_2}\in \Gamma(\cD_{-2})$,\dots,
$X_{n_{d-1}+1}, \dots, X_{n}\in \Gamma(\cD_{-d})$}\\
& \text{ such that }\; \mathop{\sum}\limits_{1\leq i\leq d}\quad
i\mathop{\sum}\limits_{n_{i-1}< j\leq n_{i}} a_j\leq k\},
\end{array}\eqno{(2.2)}
$$
where $\Gamma(\cD_{-j})$ is the space of germs at $pt$ of sections
of the bundle $\cD_{-j}$. Now, to a filtration
$$
\cD= \cD_{-1}\subset \cD_{-2} \subset \cD_{-3} \dots \subset
\cD_{-d}=TM,
$$
we assign the associated graded bundle
$$
\gr(TM)=\oplus\gr\cD_{-i},\;\text{ where
$\gr\cD_{-i}=\cD_{-i}/\cD_{-i+1}$}
$$
and the bracket of sections of $\gr(TM)$ is, by definition, the
one induced by bracketing vector fields, the sections of $TM$. We
assume a \lq\lq transitivity condition\rq\rq: The Lie algebras
$$
\gr(TM)|_{pt}\eqno{(2.3)}
$$
induced at each point $pt\in M$ are isomorphic.

The grading of the coordinates determines a
nonstandard grading of $\fvect(n)$ (recall (\ref{n_i})):
$$
\renewcommand{\arraystretch}{1.4}
\begin{array}{l}
\deg x_1=\ldots =\deg x_{n_1}=1,\\
\deg x_{n_1+1}=\ldots =\deg x_{n_2}=2,\\
\dotfill \\
\deg x_{n-n_{d-1}+1}=\ldots =\deg x_{n}=d.
\end{array}\eqno{(2.4)}
$$
Denote by $\fv=\mathop{\oplus}\limits_{i\geq -d}\fv_i$ the algebra
$\fvect(n)$ with the grading $(2.4)$. One can show that the
\lq\lq complete prolong'' of $\fg_-$ to be defined shortly, i.e.,
$(\fg_-)_*:=(\fg_-, \tilde \fg_0)_*\subset \fv$, where $\tilde
\fg_0:=\fder_0\fg_-$, preserves $\cD$.

For nonholonomic manifolds, an analog of the group $G$ from the
term ``$G$-structure'', or rather of its Lie algebra,
$\fg=\text{Lie}(G)$, is the pair $(\fg_-, \fg_0)$, where $\fg_0$
is a  subalgebra of the $\Zee$-grading preserving Lie algebra of
derivations of $\fg_-$, i.e., $\fg_0 \subset \fder_0\,\fg_-$. If
$\fg_0$ is not explicitly indicated, we assume that $\fg_0
=\fder_0\,\fg_-$, i.e., is the largest possible.

Given a pair $(\fg_-, \fg_0)$ as above, define its {\it
Tanaka-Shchepochkina prolong} to be the maximal subalgebra
$(\fg_-, \fg_0)_*=\mathop{\oplus}\limits_{k\geq -d} \fg_k$ of
$\fv$ with given non-positive part $(\fg_-, \fg_0)$. For an
explicit construction of the components, see \cite{Sh14},
\cite{Y}, \cite{ShN} and below.

\ssec{Natural bases in $T_mM$: the $D_i$'s and the $Q_i$'s
(\cite{ShN})}\index{$D_i$, vector fields} \index{$Q_i$, vector fields}
\index{vector fields, $D_i$'s and $Q_i$'s} Vershik and Gershkovich showed \cite{VG}
that every nonholonomic structure $\cD$ on $M$ determines a  structure of $\Zee$-graded
nilpotent Lie algebra in $\gr(TM)$. We will only consider  manifolds with a transitive
action of the diffeomorphism group of $M$, i.e.,
the manifolds for which these Lie algebras are isomorphic.

A natural basis in every tangent space $T_mM$ to any manifold $M$ is given by partial
derivatives. If $M$ is endowed with a nonholonomic structure,
then there are {\bf two} types of natural bases in $\gr\,T_mM$.
In physics literature on supersymmetry and supergravity, the elements of
these two bases that generate the Lie algebra $\gr\,T_mM$ are denoted the $D_i$'s
and the $Q_i$'s, respectively.

Let us consider the simplest example. Let $\dim M=2n+1$ and let
the nonholonomic structure on $M$ be given by the contact form
$\alpha=dt-\sum\left(p_idq_i-q_idp_i\right)$. The vector fields
that {\bf belong} to the distribution $\cD$ are the fields
$$
X=f\partial_t+\sum\left(g_i\partial_{q_i}+h_i\partial_{p_i}\right)\;\text{ such
that }\;
\alpha(X)=f-\sum\left(p_ig_i+q_ih_i\right)=0.\eqno{(2.5)}
$$
In particular, we see that neither $\partial_{q_i}$ nor $\partial_{p_i}$ belongs
to $\cD$, but rather
$$
D_{p_i}=\partial_{q_i}+{p_i}\partial_t\;\text{ and }\;
D_{q_i}=\partial_{p_i}-{q_i}\partial_t.
$$
These $D_{p_i}$ and $D_{q_i}$ are examples of the $D$-type basis
vectors. They, and their brackets, span the space of sections of
$gr(TM)$ at any given point $m$. By abuse of speech, we say that
the $D$-vectors span $T_mM$, and same applies to $Q$-vectors
defined below.

Now, the Lie algebra that {\bf preserves} $\cD$ consists of vector
fields $X$ such that (here $L_X$ is the Lie derivative along $X$)
$$
L_X(\alpha)=0.\eqno{(2.6)}
$$
The corresponding vector fields in our particular case of the contact
distribution are contact vector fields $K_f$ generated by $f\in\Cee[t, p, q]$:
$$
K_f=(2-E)(f)\pder{t}-H_f +
\pderf{f}{t} E, \eqno{(2.7)}
$$
where $E=\sum\limits_i y_i \pder{y_{i}}$ (here the $y_{i}$ are all
the coordinates except $t$) is the {\it Euler operator}, and $H_f$
is the Hamiltonian field with Hamiltonian $f$ that preserves
$d\alpha$:
$$
H_f=\sum\limits_{i\leq n}\left(\pderf{f}{p_i}
\pder{q_i}-\pderf{f}{q_i}
\pder{p_i}\right).
$$
(As is easy to see,
$$
L_{K_f}(\alpha)=2\displaystyle\pderf{f}{t}\alpha, \eqno{(2.8)}
$$
where $L_X$ is the Lie derivative along $X$.)
The basis of the tangent space is spanned by
$$
K_{p_i}=\partial_{q_i}-{p_i}\partial_t\;\text{ and }\;
K_{q_i}=\partial_{p_i}+{q_i}\partial_t
$$
and their brackets. These $K_{p_i}$ and $K_{q_i}$ are examples of the $Q$-type
basis vectors.

How to interpret the $D$-type and the $Q$-type vectors? Let
$$
\fn = \mathop{\oplus}\limits _{-d\leq i\leq -1}\fn_i
$$
be a nilpotent Lie algebra generated by $\fn_{-1}$. Let $B=\{b_1, \ldots, b_n\}$
be a graded basis of $\fn$ (the basis is said to be {\it graded} if its first
$n_1:=\dim \fn_{-1}$ elements span $\fg_{-1}$, the next $n_2:=\dim \fn_{-2}$
elements span $\fn_{-2}$, and so on). Let $N$ be the connected and simply
connected Lie group with the Lie algebra $\fn$. On $N$, consider the two systems
of vector fields: the left-invariant fields $D_i$ and the right-invariant fields
$Q_i$ such that ($e$ is the unit of $N$)
$$
D_i(e)=Q_i(e)=b_i\;\text{ for all }i=1, \ldots, n.
$$
{\bf NB}: Here we deviate from the conventions of physical papers where the
symbols $D_i$ and $Q_i$ are only applied to the generators of $\fn$, i.e., to
the first $n_1$ elements.

Let $\fg_-$ be a realization of $\fn$ by {\bf left-invariant} vector fields, so
the vectors $D_i(e)$ span $\fg_-$. Let $\theta^i$ be {\bf right-invariant}
1-forms on $N$ such that
$$
\theta^i(Q_j)=\delta^i_j.
$$
Now, any vector field $X$ on $N$ is of the form
$$
X=\sum_{i=1}^n \theta^i(X) Q_i.\eqno{(2.9)}
$$
Since each $D_i$ commutes with each $Q_j$ (if $\fn$ is a Lie superalgebra, they
{\bf super}commute), it follows that
$$
\theta^i([D_j, X])=D_j(\theta^i(X)).
$$

Now, let us determine a {\bf right-invariant} distribution $\cD$
on $N$ such that $\cD|_e=\fn_{-1}$. Clearly, $\cD$ is singled out
in $TN$ by eqs. for $X\in \fvect(n)$
$$
\theta^{n_{1}+1}(X)=0, \quad \dots, \quad \theta^{n}(X)=0.
$$

Since each $D_i$ commutes with each $Q_j$, the algebra $\fg_-$ preserves $\cD$.
The coordinates $(2.4)$ on $N$ described above determine two embeddings of $\fn$
into $\fvect(n)$: one is spanned by the $D_i$ and the other one by the $Q_i$.

Denote by $\fg=\mathop{\oplus}\limits_{i\geq -d}\fg_i$ the algebra
$\fvect(n)$ with the grading $(\ref{gr})$. Then
$\fg_-=\mathop{\oplus}\limits_{i<0}\fg_i$ preserves $\cD$. We will
show later that the \lq\lq complete prolongation'' of $\fg_-$,
i.e., $(\fg_-)_*:=(\fg_-, \tilde \fg_0)_*$, where $\tilde
\fg_0:=\fder_0\fg_-$, also preserves $\cD$.

Thus we see that, with every nonholonomic manifold $(M, \cD)$, a
natural $G$-structure is associated, its Lie algebra is
$\Lie(G)=\fder_0\,\fg_-$. But the structure functions of this
$G$-structure do not reflect the nonholonomic nature of $M$.

Indeed, recall an example from \cite{St}. Let $W_1\subset W$ be a
subspace of dimension $k$ and $G\subset GL(W)$ the parabolic
subgroup that preserves the subspace. Then to determine a
$G$-structure on $M$, where $\dim M=\dim W$, is the same as to
determine a differential $k$-system or a $k$-dimensional
distribution. A fixed frame $f$ in $T_mM$ determines an
isomorphism $f: W\tto T_mM$. Given a $G$-structure on $M$, we set
$\cD(m)=f(W_1)$. Since $G$ preserves $W_1$, the subspace $\cD(m)$
indeed depends only on $m$, not on $f$.

The other way round, given a distribution $\cD$, consider the frames $f$ such
that
$f^{-1}(\cD(m))=W_1$. They form a $G$-structure. The flat $G$-structures
correspond
to integrable distributions.

To take the nonholonomic nature of $M$ into account, we need something new
--- an analog of the above Proposition 1.1 for the case where the natural basis
of the tangent space consists not of partial derivatives but
rather of covariant derivatives corresponding to the connection
determined by the same Pfaff equations that determine the
distribution, and therefore instead of $T_{m}M=\fg_{-1}$ we have
$(gr(TM))_{m}=\fg_{-}$. To be able to formulate such Proposition,
we have to define
\medskip
(1) the simplest nonholonomic structure --- the ``flat'' one,

(2) the analog of $\fg_{0}$ when $\fg_{-1}$ is replaced by
$\fg_{-}$ and only distribution is given,

(3) what is the analog of $(\fg_{-1}, \fg_0)_*$,

(4) what is the analog of $H^{k, 2}_{(\fg_{-1}, \fg_0)_*}$.

\medskip
Here are the answers:

1) Let $\cD$ be a nonholonomic distribution in $M$, let $F$ be the
flag which $\cD$ determines at a point $m\in M$. Let $N:=\Kee^n$
with a fixed flag $F$ and a fixed frame $f$. Having identified
$T_nN$ with $N$ by means of the translation by $n$ considered as
an element of the {\it nilpotent} Lie group $N$ whose Lie algebra
is $\fg_-$ (since the group $N$ is not commutative now, we select,
say, left translations) we fix a frame and a flag --- the images
of $f$ and $F$ --- in each $T_nN$. A {\it flat nonholonomic
structure} on $N$ is the {\it pair} of bundles (the frame bundle,
the distribution $\cD$); the fibers of both bundles over $n$ are
obtained from the fixed frame and flag, respectively, by means of
the $G$-action, where $G$ is the (connected and simply connected)
Lie group whose Lie algebra $\fg_0$ is defined at the next step.

2) If only a distribution $\cD$ is given, we set $\fg_0:=\fder_0\fg_-$; it is
often interesting to consider an additional structure on the distribution,
say Riemannian, cf. \cite{VG2}, as in the case of Carnot-Carath\'eodory metric
in which case $\fg_0$ is a subalgebra of $\fder_0\fg_-$, e.g., $\fder_0\fg_-\cap \fo(\fg_{-1})$.

3) Given a pair $(\fg_-, \fg_0)$ as above, define its $k$th {\it
Tanaka-Shchepochkina prolong} (given simultaneously, although
\cite{Sh1} was published later than \cite{T}; \cite{Sh1} also
embraces Lie superalgebras and various {\it partial prolongs}, see
\cite{Sh2}) for $k>0$ to be:
$$
\fg_k = (i(S^{\bcdot}(\fg_-^*)\otimes \fg_0)\cap j(S^{\bcdot}(\fg_-^*)\otimes
\fg_-))_k,
\eqno{(2.10)}
$$
where the subscript singles out the component of degree $k$, where
$S^{\bcdot}=\oplus S^i$ and $S^i$ denotes the operator of the
$i$th symmetric power, and where
$$
\renewcommand{\arraystretch}{1.4}
\begin{array}{l}
 i: S^{k+1}(\fg_{-1}^*)\otimes \fg_{-1}\tto S^{k}(\fg_{-1}^*)\otimes
\fg_{-1}^*\otimes\fg_{-1}, \\
j:
S^{k}(\fg_{-1}^*)\otimes \fg_{0}\tto S^{k}(\fg_{-1}^*)\otimes
\fg_{-1}^*\otimes\fg_{-1}
\end{array}
$$
are natural embeddings.

Similarly to the case where $\fg_-$ is commutative, define $(\fg_-,
\fg_0)_*$ to be $\mathop{\oplus}\limits_{k\geq -d} \fg_k$ with $\fg_k$ for
$k>0$ given by $(2.10)$; then, as is
easy to verify, $(\fg_-, \fg_0)_*$ is a Lie algebra.

4) Arguments similar to those of \cite{St} should show
that $H^2(\fg_-; (\fg_-, \fg_0)_*)$ is the space of
values of all nonholonomic structure functions --- obstructions to the
identification of the infinitesimal neighborhood of a point $m$ of the
manifold $M$ with a nonholonomic structure (given by $\fg_-$ and
$\fg_0$) with the infinitesimal neighborhood of a point of a flat
nonholonomic manifold with the same $\fg_-$ and $\fg_0$.
We intend to give a detailed proof of this statement elsewhere.

The space $H^2(\fg_-; (\fg_-, \fg_0)_*)$ naturally splits into
homogeneous components whose degrees will be called the {\it orders}
of the structure functions; the orders run $2-d$ to $N+2d$ (or to
$\infty$ if $N=\infty$). As in the case of a commutative
$\fg_-=\fg_{-1}$, the structure functions of order $k$ can be
interpreted as obstructions to flatness of the nonholonomic manifold
with the $(\fg_-, \fg_0)$-structure provided the obstructions of
lesser orders vanish. Observe that, for nonholonomic manifolds, the
order of structure functions is no more in direct relation with the
orders of the infinitesimal neighborhoods of the points we wish to
identify: distinct partial derivatives bear different \lq\lq degrees''.

Different filtered algebras $L$ with the same graded $\fg_-$ are
governed precisely by the coboundaries responsible for filtered
deformations of $\fg_-$, and all of them vanish in cohomology, so the above
nonholonomic structure functions are well-defined.

\section*{\S 3. The Riemann and Weyl tensors. Projective structures}

\ssec{The conformal case} For the classical domains $X=S/P$ that Goncharov
considered, the structure functions are generalizations of the Weyl
tensor --- the conformally invariant part of the Riemann tensor (the
case $S=O(n+2)$ and $G=CO(n)$). In most of these cases
$$
(\fg_{-1}, \fg_0)_*=\fs\; (:=\text{Lie}(S))\eqno{(3.1)}
$$
and the description of the structure functions is a particular case of
the BWB theorem. In particular, if $(3.1)$ holds, the space
$H^2(\fg_{-1}; (\fg_{-1}, \fg_0)_*)$, considered as a $\fg_0$-module,
has the same number of irreducible components and the same dimension
as $E^2(\fg_{-1})$; only weights differ.

\ssec{The generalized Riemannian case} When we reduce $\fg_{0}$,
by retaining its semi-simple part $\hat{\fg_{0}}$ and deleting the
center, we can not directly apply the BWB theorem because
$(\fg_{-1}, \hat{\fg_{0}})_*=\fg_{-1}\oplus \hat{\fg_{0}}$ is not
simple but we can reduce teh problem to the conformal case, since,
as is known,
$$
H^2(\fg_{-1}; (\fg_{-1}, \hat{\fg_{0}})_*)= H^2(\fg_{-1}; \fs)\oplus
S^2(\fg_{-1}^*).\eqno{(3.2)}
$$

For the nonholonomic case, a similar reduction is given by
Premet's theorem (below). Its general case, though sufficiently
neat, is not as simple as (3.2). However, although the following
analog of (3.2) is not always true
$$
H^2(\fg_{-}; (\fg_{-}, \hat{\fg_{0}})_*)= H^2(\fg_{-}; \fs)\oplus
S^2(\fg_{-1}^*),\eqno{(3.3)}
$$
it is still true in many cases of interest: for the \lq\lq contact
grading\rq\rq.

\ssec{The projective case. Theorems of Serre and Yamaguchi} When
$(3.1)$ fails, $\fs$ is a proper
subalgebra of $(\fg_{-1}, \fg_0)_*$. It is of interest therefore

(a) to list all the
cases where, having started from a simple Lie (super)algebra
$\fs=\mathop{\oplus}\limits_{i\geq -d}\fs_i$, we have the following
analog of (3.1)
$$
(\fs_-, \fs_0)_*=\fs\eqno{(3.4)}
$$ and

(b) find out {\it what} is the \lq\lq complete prolongation'' of
$\fs_-$, i.e., what is $(\fs_-)_*:=(\fs_-, \tilde \fs_0)_*$, where
$\tilde \fs_0:=\fder_0\fs_-$.

For the simple finite dimensional Lie algebras $\fs$, Yamaguchi
\cite{Y} gives the answer. It is rather interesting and we
reproduce it. The answer for simple Lie superalgebras is obtained
by Shchepochkina (unpublished). Comment: one would expect that
$\tilde \fs_0$ strictly contains $\fs_0$, and hence $(\fs_-)_*$
should strictly contain $\fs$; instead they are equal (in
particular, $\tilde \fs_0=\fs_0$).

\begin{Theorem}[\cite{Y}] Equality $(\fs_-)_*=\fs$ holds almost always. The
exceptions are

{\em 1)} $\fs$ with the grading of depth $d=1$ (in which case
$(\fs_-)_*=\fvect(\fs_-^*)$);

{\em 2)} $\fs$ with the grading of depth $d=2$ and $\dim\fs_{-2}=1$,
i.e., with
the \lq\lq contact'' grading, in which case $(\fs_-)_*=\fk(\fs_-^*)$ (these
cases
correspond to exclusion of the nodes on the Dynkin graph connected with the node
for the maximal root on the extended graph);

{\em 3)} $\fs$ is either $\fsl(n+1)$ or $\fsp(2n)$ with the
grading determined by \lq\lq selecting'' the first and the $i$th
of simple coroots, where $1<i<n$ for $\fsl(n+1)$ and $i=n$ for
$\fsp(2n)$. (Observe that, in this case, $d=2$ with $\dim
\fs_{-2}>1$ for $\fsl(n+1)$ and $d=3$ for $\fsp(2n)$.)

Moreover, $(\fs_-, \fs_0)_*=\fs$ is true almost always. The cases
where this fails (the ones where a projective action is possible)
are $\fsl(n+1)$ or $\fsp(2n)$ with the grading determined by
\lq\lq selecting'' only one (the first) simple coroot.
\end{Theorem}

\underline{Case 1) of Yamaguchi's theorem}: for the conformal (Weyl)
case, see \cite{Go}; for the Riemannian case, see \cite{LPS} .

For the classical domains $X=S/P$, (3.1) fails only for $S=SL(n+1)$ and
$X=\Cee P^n$;
in this case $\fg_{0}=\fgl(n)$ and $(\fg_{-1}, \fg_0)_*=\fvect(n)$, the
Lie algebra
of vector fields in $n$ indeterminates. The space of \lq\lq total''
structure functions
$H^2(\fg_{-1}; \fvect(n))$ differs from $H^2(\fg_{-1}; \fs)$, the latter
structure
functions correspond to obstructions to the {\it projective structure}.
For many
facets of projective structures, see \cite{OT} and \cite{BR}.

The Riemannian version of this projective case, corresponds to
$\hat{\fg_{0}}=\fsl(n)$ and $(\fg_{-1}, \hat{\fg_{0}})_*=\fsvect(n)$, the
Lie algebra of
divergence free vector fields.

The cases of \lq\lq complete prolongation''
$(\fs_{-1})_*=\fvect(\fs_{-1}^*)$ and their \lq\lq Riemannian
version'' $(\fs_{-1})_*=\fsvect(\fs_{-1}^*)$, as well as
$(\fs_{-1})_*=\fh(\fs_{-1}^*)$, were considered by Serre long ago,
see \cite{St}, and the answer is as follows:

\begin{Theorem}[Serre, see \cite{St}; for super version, see \cite{LPS} and
\cite{GLS}] {}~{}

{\em 1)} $H^2(\fs_{-1};
\fvect(n))=0$ and $H^2(\fs_{-1};
\fsvect(n))=0$.

{\em 2)} $H^2(\fs_{-1}; \fh(2n))=E^3(\fs_{-1}^*)$. \end{Theorem}

\underline{Case 2) of Yamaguchi's theorem} is taken care of by one of Premet's
theorems and formula (3.5) below.

\underline{Case 3) of Yamaguchi's theorem} is done in \S 6 of this paper.

In what follows, for manifolds $X=S/P$ with nonholonomic
structure, we say ``N-Weyl'' or ``N-conformal'', for tensors
corresponding to cohomology of $\fg_{-}$ with coefficients in
$(\fg_{-}, \fg_0)_*$, ``N-Riemannian'' for nonholonomic structure
functions $(\fg_{-}, \hat{\fg_{0}})_*$, where $\hat{\fg_{0}}$ is
the semi-simple part of $\fg_{0}$, and ``N-projective'' for the
coefficients in $\fs=\text{Lie}(S)$ whenever $\fs$ is smaller than
$(\fg_{-}, \fg_0)_*$, for example, for {\it partial Cartan
prolongs}, see \cite{LSh}.

\ssec{The simplest examples (exclusion of the first simple coroot of
$\fsp(2n+2)$)} Let $\fg_-=\fhei(2n)$, the Heisenberg Lie
algebra. Then $\fg_{0}=\fc\fsp(2n)$ (i.e., $\fsp(2n)\oplus \Cee z$) and
$(\fg_{-},
\fg_{0})_*$ is the Lie algebra $\fk(2n+1)$ of contact vector fields.

So far, there is no analog of Serre's theorem on involutivity for
simple $\Zee$-graded Lie algebras of depth $>1$, cf. \cite{LPS}, and
examples from \cite{GLS} show that if exists, the theorem is much more
involved.

The fact that
$$
H^2(\fhei(2n); \fk(2n+1))=0 \eqno{(3.5)}
$$
explains why {\bf the Pfaff equation $\alpha(X)=0$ for
$X\in\fvect(2n+1)$ can be reduced to a canonical form}, cf.
\cite{Z}. This fact is an easy corollary of a statement on
cohomology of coinduced modules \cite{FF}. For the N-Riemannian
tensor in this case, we have: $\hat{\fg_{0}}=\fsp(2n)$ and
$(\fg_{-}, \hat{\fg_{0}})_*$ is the Poisson Lie algebra
$\fpo(2n)$. The Poisson Lie algebra is spanned by fields $K_f$,
where $\displaystyle\pderf{f}{t}=0$. Now, from $(3.5)$ and the
short exact sequence
$$
0\tto \fpo(2n)\tto \fk(2n+1)\stackrel{\pder{t}:K_f\mapsto
\pderf{f}{t}}{\tto}\Cee[t, p, q]\tto 0
$$
we easily deduce (using the corresponding long exact sequence, see
\cite{FF}) that
$$
H^2(\fhei(2n); \fpo(2n))=0. $$ In our terms, this fact (usually
called {\it Darboux's theorem} and proved by analytic means
\cite{Z}) is an explanation why {\bf the contact form $\alpha$ can
be reduced to a canonical form not only at any point but locally}.

\ssec{Other examples} For numerous examples of N-projective structures in
various instances, see
 \cite{C1}--\cite{C5} and \cite{YY}, and (in super setting) \cite{MaG}. Armed
with {\bf SuperLie}, one can now easily perform the computations
of relevant Lie algebra cohomology. Premet's theorems tell what to
compute in the N-Riemannian case and again with {\bf SuperLie}
this will be easy: we just give a few samples ({\bf one} selected
simple coroot for every $\fs$ and {\bf two} selected coroots for
the two series of one of Yamaguchi's cases).

\section*{\S 4. Premet's Theorems (from Premet's letter to DL, 10/17/1990)}

In 1990, DL asked Alexander Premet: how to reduce computations of the space of
values for a nonholonomic Riemann tensors to that for the
nonholonomic Weyl tensor, as in (3.2)? Namely, is $(3.3)$ always true?

Premet wrote two letters with a general answer. One letter is
reproduced practically without changes below (DL is responsible for any mistakes
left/inserted); it shows how to reduce the problem to
computing (the 1st) cohomology of $\fg_{-}$ with coefficients in a
certain $\fg_{-}$-module which is not a $\fg$-module. Little was
known about such cohomology except theorems of Kostant (on $H^1$) and
of Leger and Luks (on $H^2$) both for the case where $\fg_{-}$ is the maximal
nilpotent subalgebra. Premet's second letter (reproduced in \cite{LLS})
contained a mighty generalization of these theorems for $H^i$ for any
$i$ and any $\fg_{-}$.

However, in nonholonomic cases, to derive an {\it explicit} answer
from the BWB theorem is difficult ``by hands'', the extra terms in
the Riemannian case (see sec.~ 4.4 below) add extra job. So
Premet's theorems were put aside for 13 years. Now that a package
{\bf SuperLie} (\cite{Gr}), originally designed for the purposes
of supergravity, is sufficiently developed, we are able to give an
explicit answer: see the next section. The cases we consider here
(of the maximal parabolic subalgebras) required several minutes to
compute. (But much longer to document the results, and it will require
 a while to interpret them, say as in \cite{LPS}.)
To our regret, Premet looks at his theorems as a mere
technical exercise (\lq\lq a simple job for Kostant'') not
interesting enough to co-author the paper.

\ssec{4.1. Terminological conventions} Let $\fg$ be a simple (finite
dimensional) Lie algebra. Let $L^\lambda$ denote the irreducible
(finite dimensional) $\fg$-module with the highest weight $\lambda$;
let $E_{\mu}$ be the subspace the module $E$ of weight $\mu$.

Let $R$ be the root system of $\fg$ and $B$ the base (system of simple
roots). Let $W=W(R)$ be the Weyl group of $\fg$ and $l(w)$ the
length of the element $w\in W$; let $W_{i}$ be the subset of elements
of length $i$. Let $R_I\subset R$ and let $B_I$ be the base of $R_I$. Set
(this is a definition of $k(i)$ as well)
$$
W(I)_{i}=\{w_{i, 1}, \ldots , w_{i, k(i)}\in W_{i}\mid w_{i,
j}^{-1}(B\setminus B_{I})>0\;\text{ for all $1\leq j\leq
k(i)$}\}.\eqno{(4.1)}
$$

Let the Dynkin graph of $B$ be, for example, as follows:
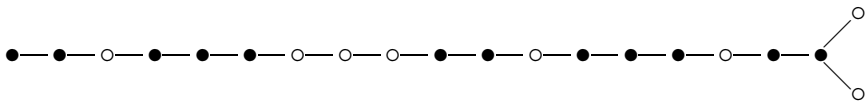
\begin{center}
\begin{picture}(400, 50)
\put(0, 25){$\bullet$}
\put(6, 28){\line(1, 0){10}}
\put(18, 25){$\bullet$}
\put(24, 28){\line(1, 0){10}}
\put(36, 25){$\circ$}
\put(42, 28){\line(1, 0){10}}
\put(54, 25){$\bullet$}
\put(60, 28){\line(1, 0){10}}
\put(72, 25){$\bullet$}
\put(78, 28){\line(1, 0){10}}
\put(90, 25){$\bullet$}
\put(96, 28){\line(1, 0){10}}
\put(108, 25){$\circ$}
\put(114, 28){\line(1, 0){10}}
\put(126, 25){$\circ$}
\put(132, 28){\line(1, 0){10}}
\put(144, 25){$\circ$}
\put(150, 28){\line(1, 0){10}}
\put(162, 25){$\bullet$}
\put(168, 28){\line(1, 0){10}}
\put(180, 25){$\bullet$}
\put(186, 28){\line(1, 0){10}}
\put(198, 25){$\circ$}
\put(204, 28){\line(1, 0){10}}
\put(216, 25){$\bullet$}
\put(222, 28){\line(1, 0){10}}
\put(234, 25){$\bullet$}
\put(240, 28){\line(1, 0){10}}
\put(252, 25){$\bullet$}
\put(258, 28){\line(1, 0){10}}
\put(270, 25){$\circ$}
\put(276, 28){\line(1, 0){10}}
\put(288, 25){$\bullet$}
\put(294, 28){\line(1, 0){10}}
\put(306, 25){$\bullet$}
\put(310, 31){\line(1, 1){10}}
\put(310, 25){\line(1, -1){10}}
\put(320, 41){$\circ$}
\put(320, 10){$\circ$}
\end{picture}
\end{center}
\noindent and let $B_I$ consist of roots corresponding to the black
nodes. Let us represent $B_I$ as the union of connected subgraphs:
$$
B_I = B_I^{(1)}\coprod \ldots \coprod B_I^{(s)},
$$
where $s$ (in our example $s=5$) is the number of connected
components of the Dynkin graph $D_I$ of $B_I$ and where
$B_I^{(i)}$ corresponds to the $i$th connected components of $D_I$
(counted from left to right). Set
$$
c = \card~B, \quad c_i = \card~\{\alpha \in B\setminus B_I\mid (\alpha ,
B_I^{(i)}) \neq 0\} - 1.
$$
Clearly, if $B_I \not = B$, then $c_i \in \{0, 1, 2\}$. For example,
for the graph of $\fo(20)$ depicted above, we have:
$$
c = 20, \; \; c_1 = 0, \; \; c_2 = 1, \; \; c_3 = 1, \; \; c_4 = 1, \; \;
c_5 = 2.
$$
The following statement is obvious.

\begin{Statement} {\em 1)} $c_i = 2$ if and only if $R$ is of type
$D_n$, $E_6$, $E_7$, $E_8$, one of the endpoints of $D_I^{(i)}$ is
a branching point for $D$, and the remaining endpoint of $D_I^{(i)}$
is not an endpoint of $D$.

{\em 2)} $c_i = 0$ if and only if all but one of the end vertices of
the graph of $B_I^{(i)}$ are the end vertices for the graph of
$R$.
\end{Statement}

\ssec{4.2. The Borel-Weil-Bott theorem} Let $\rk \fg =r>1$,
$I\subset \{1, \dots, r\}$; let $\fp =\fp_{I}$ be a parabolic
subalgebra generated by the Chevalley generators $X_{i}^{\pm}$ of
$\fg$ except the $X_{i}^{+}$, where $i\in I$. As is known, $\fp
=\fg_-\oplus \fl$, where $\fl$ is the Levi (semi-simple)
subalgebra generated by all the $X_{i}^{\pm}$, where $i\not\in I$.
Clearly, $\fl=\fl^{(1)}\oplus\fz$, where $\fl^{(1)}$ is the
derived algebra of $\fl$, and $\fz = \fz (\fl)$ is the center of
$\fl$.

So, in terms of \S 3, $\fg_{0}=\fl$, $\hat\fg_{0}=\fl^{(1)}$.

\begin{Theorem}[The BWB Theorem, see \cite{BGG}] Let $E = L^{\lambda}$
be an irreducible (finite dimensional) $\fg$-module with {\bf highest}
weight $\lambda$. Then $H^i(\fg_-; E)$ is the direct sum of
$\fl$-modules with the {\bf lowest} weights $-w_{i, j}(\lambda
+\rho)+\rho$, where $w_{i, j}\in W(I)_{i}$, see $(4.1)$; each such module
enters with multiplicity $1$.
\end{Theorem}

The BWB theorem describes (for $i=2$) nonholonomic analogs of Weyl
tensors. Theorem 4.4 describes nonholonomic analogs of the Riemann
tensors. To prove sec.~4.4, we need the following Lemma.

\ssbegin{4.3}{Lemma} Let $E = L^{\lambda}$ be an irreducible (finite
dimensional) $\fg$-module with highest weight $\lambda$ such that $E
\simeq E^*$. Let $V$ be a $\fp$-invariant subspace in $E$ which
contains $E_- = \mathop{\oplus}\limits_{\mu= \sum k_i\alpha _i\mid k_i
< 0}E_{\mu}$. Let further any weight $\mu$ of
$E/V$ be of the form $\mu=-\sum a_i\alpha_i$, where $a_i\ge 0$.

Then, for any $i < \rk \fg$, we have the $\fl$-module
isomorphism:
$$
H^i(\fg_-; V) \simeq H^i(\fg_-; E)\oplus H^{i-1}(\fg_- ; E/V).
$$
\end{Lemma}

\begin{proof} As is well known \cite{FF}, with every short exact sequence
of
$\fg_-$-modules
$$
0 \tto V \tto E \tto E/V \tto 0
$$
there is associated the long exact cohomology sequence
$$
\renewcommand{\arraystretch}{1.4}
\begin{array}{l}
    0 \tto H^0(\fg_-; V) \tto H^0(\fg_-; E)
\buildrel\varphi _0\over \tto H^0(\fg_-; E/V) \tto \\
H^1(\fg_-; V)
\tto H^1(\fg_-; E)
\buildrel\varphi _1\over \tto H^1(\fg_-; E/V) \tto \ldots\\
\tto H^i(\fg_-; V) \tto H^i(\fg_-; E)
\buildrel\varphi _i\over \tto H^i(\fg_-; E/V)
\tto \ldots
\end{array}
$$
Let us prove that the weight $-w(\lambda+\rho)+\rho$ can not be a
weight of the $\fl$-module $H^i(\fg_-; E/V)$ for $i <\rk~\fg$.
Indeed, each weight of $H^i(\fg_-; E/V)$, which is a quotient of a
submodule of $E^i(\fg_-^*)\otimes E/V$, is of the form
$$
\gamma _1 + \ldots + \gamma _k + \mu,
$$
where $\gamma _1, \ldots , \gamma _k$ are distinct positive roots
which do not belong to the root lattice $Q(B_I)$, and $\mu $ is a
weight of $E$ of the form $\mu = \mathop{\sum}\limits_{a_i\geq 0}
a_i\alpha _i$.

Suppose that
$$
-w(\lambda+\rho)+\rho = \gamma _1 + \ldots + \gamma _k + \mu.
$$
Then $-\mu -w(\lambda)-w(\rho)+\rho = \gamma _1 + \ldots + \gamma _k
$. Set
$$
R_{W}^{-} = \{a\in R_+\mid w^{-1}(\alpha) > 0\}\text{ and }
R_{W}^{+} = \{a\in R_+\mid w^{-1}(\alpha) < 0\}.
$$
Let $\gamma
_1, \ldots , \gamma _s \in R_{W}^{-}$ and $\gamma _{s+1}, \ldots , \gamma _k
\in R_{W}^{+}$.

Recall that
$$
\rho - w(\rho ) =\mathop{\sum}\limits_{\beta \in R_{W}^{-}}\beta .
$$
Since $E = E^*$, the weight $-\mu $ is a weight of $E$; but then
$$
\renewcommand{\arraystretch}{1.4}
\begin{array}{l}
w^{-1}(\gamma _1 + \ldots + \gamma _k) = \\
\left(w^{-1}(\gamma _1) + \ldots +
w^{-1}(\gamma _s)\right) + \left(w^{-1}(\gamma _{s+1}) + \ldots + w^{-1}
(\gamma
_k)\right) =\\
(-\lambda - w^{-1}(\mu )) + w^{-1} (\mathop{\sum}\limits_{\beta \in
R_{W}^{-}}\beta).
\end{array}
$$
Since $-w^{-1}(\mu )$ is a weight, $\lambda + w^{-1}(\mu )\in Q_+(R)$.
Hence,
$$
-(\lambda + w^{-1}(\mu )) + w^{-1} (\mathop{\sum}\limits_{\beta \in
R_{W}^{-}}\beta ) \in -Q_+(R).
$$
On the other hand,
$$
w^{-1}(\gamma _{s+1}) + \ldots + w^{-1}(\gamma _k)\in Q_+(R)
$$
and all the $w^{-1} (\gamma _i)$ for $i\leq s$ enter, as summands,
$\mathop{\sum}\limits_{\beta \in R_{W}^{-}}w^{-1}\beta$, and
therefore cancel each other.

We finally get:
\begin{multline*}
Q_+(R)\ni w^{-1}(\gamma _{s+1}) + \ldots + w^{-1}(\gamma _{k}) =\\
-(\lambda + w^{-1}(\mu)) + w^{-1}(\mathop{\sum}\limits_{\beta \in
R_{W}^{-} \setminus \{\gamma _1, \ldots , \gamma _s\}}\beta ) \in
-Q_+(R).
\end{multline*}

Thus, $s = k, \lambda = -w^{-1} (\mu )$. In other words,
$$
-\mu = w(\lambda ) = -\sum a_i\alpha _i\; \text{ for }\; a_i \geq 0.
$$
Since $\lambda $ is a dominant weight, $\lambda = \sum m_i\alpha
_i$, where all the $m_i\in \Qee$ are positive. This is true for any
fundamental weight, as follows from the tables from \cite{Bu}.

By the hypothesis, $l(w) < \rk ~\fg$, and therefore
$$
w = s_{\alpha _{1}}\ldots s_{\alpha _{r}}, \; \text{ where }\; \; r <
\rk ~\fg.
$$
Set
$$
B_{w}:=\{\alpha _{i_{1}}, \ldots , \alpha
_{i_{r}}\}\subset B.
$$
For any $x \in Q(B)\otimes_{\Zee}\Ree$, we have $x - w(x) =
-\mathop{\sum}\limits_{\alpha \in B_{w}} n_{\alpha}\alpha$. But
then (recall that $\varpi_j$ is the $j$th fundamental weight, see
\cite{Bu})
$$
(x - w(x), \varpi _j) = 0\; \; \text{ for some } \; j\leq \rk ~\fg.
$$
Therefore, $(\lambda , \varpi _j) = (w(\lambda), \varpi _j)$. Now
notice that
$$
(\lambda, \varpi _j) = m_j(\alpha _j, \varpi _j) > 0, \; \;
(w(\lambda) , \varpi _j) = -a_j(\alpha _j, \varpi _j) \leq 0.
$$
This contradiction shows that $\varphi _i=0$ for $i < \rk ~\fg$. This,
in turn, means that, for $i = 1, \ldots , \rk\fg-1$, there exist short
exact sequences of
$\fl$-homomorphisms
$$
0 \tto H^{i-1}(\fg_-; E/V) \tto H^i(\fg_-; V) \tto H^i(\fg_-; E) \tto
0. \qed
$$
\noqed \end{proof}

\ssbegin{4.4}{Theorem} $H^2(\fg_-; \fg_-\oplus \fl^{(1)}) =
H^2(\fg_-; \fg)\oplus H^1(\fg_-; (\fg_-\oplus \fz)^*)$.
\end{Theorem}

\begin{proof} Set $E = \fg, V = \fg_-\oplus\fl^{(1)}$. By Lemma 4.3,
$$
H^2(\fg_-; \fg_-\oplus\fl^{(1)}) = H^2(\fg_-; \fg)\oplus H^1(\fg_-;
\fg/(\fg_-\oplus\fl^{(1)})).
$$
It remains to verify that $\fg/(\fg_-\oplus\fl^{(1)}) =
(\fg_-\oplus\fz)^*$. Indeed, the Killing form $K$ establishes an
isomorphism $\fg = \fg ^*$, and therefore
$$
(\fg/(\fg_-\oplus\fl^{(1)}))^* = \{x\in \fg\mid
K(x, \fg_-\oplus\fl^{(1)}) = 0\}
= \fg_-\oplus\fz,
$$
where $\fz = \{z\in \fl\mid K(z, h_{\alpha }) = 0$ for any $\alpha \in
B_I\}$,
as required. \end{proof}

Observe that $\dim ~\fz$ is equal to the cardinality of $I$, it is 1
in \S\S 5, 7 and 2 in \S 6.

\begin{Corollary} Let $B_1 = B\setminus B_I$; let $R_1$ be the root
system generated by $B_1$ and
$$
\fg_-^{ab} = (\fg_-/\fg_-^{(1)})^* = H^1(\fg_-).
$$

{\em 1)} The following sequence is exact:
$$
\renewcommand{\arraystretch}{1.4}
\begin{array}{l}
0 \tto \fg_-^{ab} \tto \fz^*\otimes\fg_-^{ab} \tto H^1(\fg_-;
\fg/(\fg_-\oplus \fl ^{(1)})) \tto H^1(\fg_-; \fg_-^*) \tto \\
H^2(\fg_-)\oplus\mathop{\oplus}\limits_{w\in W(R_1)_{(2)}} L^{\rho-
w(\rho)} \tto 0.
\end{array}
$$

{\em 2)} If $\dim ~\fz= 1$, then the sequence
$$
0 \tto H^1(\fg_-; (\fg_-\oplus \fz)^*) \tto
H^1(\fg_-; \fg_-^*) \tto H^2(\fg_-) \tto 0
$$
is exact. In particular, if $\fg_-$ is a Heisenberg algebra (the case of contact
grading), then
$$
\text{\fbox{$H^2(\fg_-; \fg_-\oplus \fl^{(1)})\simeq H^2(\fg_-; \fg)\oplus
S^2(\fg_-/\fz(\fg_-))^*=H^2(\fg_-; \fg)\oplus S^2(\fg_{-1}^*).$}}
$$

{\em 3)} if $\fg_- =\fg_{-1}$ (is abelian), then
$$
H^2(\fg_-; \fg_- \oplus \fl^{(1)}) \simeq H^2(\fg_-; \fg)\oplus
S^2(\fg_-^*)=H^2(\fg_-; \fg)\oplus S^2(\fg_{-1}^*).
$$
\end{Corollary}

\begin{proof} From the short exact sequence
$$
0 \tto (\fg_-\oplus\fl)/(\fg_-\oplus\fl^{(1)}) \tto
\fg/(\fg_-\oplus\fl^{(1)})) \tto \fg/(\fg_-\oplus\fl)\tto 0,
$$
where $\fg/(\fg_-\oplus\fl) = \fg_-^*$ and
$(\fg_-\oplus\fl)/(\fg_-\oplus\fl^{(1)}) = \fz ^*$, we get the long
exact sequence
$$
\renewcommand{\arraystretch}{1.4}
\begin{array}{l}
    0 \tto \fg_- ^{ab} \tto \fz ^*\otimes H^1(\fg_-)
\tto H^1(\fg_-; \fg/(\fg_-\oplus\fl^{(1)})) \tto
H^1(\fg_-; \fg_-^*) \buildrel\psi\over\tto \\
     \tto   H^2(\fg_-; \fz ^*) ~~(= H^2(\fg_-)\otimes \fz^*))
\tto \ldots
\end{array}
$$
To compute the image of the map $\psi$, consider the $\fl$-module
$$
\renewcommand{\arraystretch}{1.4}
\begin{array}{ll}
M = &\{f\in \Hom ~(\fg_-, \fg/(\fg_-\oplus\fl^{(1)}))\mid Xf(Y) - Yf(X)
- f([X, Y])\in \fz^*\\
&\text{ for any }\; X, Y\in \fg_-\}.
\end{array}
$$
From the general properties of the long exact sequences (\cite{FF}) we
deduce that
$$
\IM~\psi = M/Z^1(\fg_-; \;\fg/(\fg_-\oplus\fl^{(1)})).
$$
By the BWB theorem,
$$
H^2(\fg_-) = \mathop{\oplus}\limits_{\{w\;\mid \;l(w) = 2, \; w(B\setminus
B_I)> 0\}}
L^{w(\rho ) - \rho}.
$$
Since $M$ is the direct sum of its weight subspaces relative to
$\fh$, let us study the subspaces $M_{\rho -w(\rho)}$. The weight of
$f\in M$ is equal to $w(\rho ) -\rho $ if and only if $f$ sends
$\fn_{-\gamma}$ to
$$
 (\fg/(\fg_-\oplus\fl^{(1)}))_{\rho-w(\rho)-\gamma}.
$$
Therefore, either $\rho - w(\rho ) - \gamma = \gamma '\in R$ or $\rho
- w(\rho ) - \gamma = 0$.

The second option is ruled out since $\rho - w(\rho )$ is not a root:
indeed, $l(w)=2\not =1$.

In the first case, $\rho - w(\rho ) = \gamma + \gamma '$, where
$\gamma, \gamma '\in R_+$. If $\gamma = \gamma '$, then $\rho -
w(\rho ) = 2\gamma $, implying $w^{-1} (\gamma ) < 0$. But
$\card~(R^{-}_W) = 2$ and, $\rho - w(\rho) = 2\gamma = \gamma _1 +
\gamma_2$, where $\gamma _1, \gamma_2\in R^{-}_W$. Since one of
the $\gamma_i$ is equal to $\gamma$, so is the other one. This
contradicts the hypothesis: $\gamma_1 \neq \gamma_2$.

Thus, $\rho - w(\rho ) =\gamma + \gamma '$, where $\gamma \neq \gamma
'$. It is not difficult to deduce from this that $w^{-1}(\gamma ) <
0$ and $w^{-1}(\gamma ') < 0$. Hence, $f(\fn_{-\delta }) \not = 0$ only
if $w^{-1} (\delta ) < 0$ for $f\in M_{\rho -w(\rho )}$. In this case
$f(\fn_{-\delta }) \in (\fg/(\fg_-\oplus\fl^{(1)}))_{\delta '}$,
where $w^{-1}(\delta ') < 0$. Besides, since $\delta , \delta '\not
\in R$, we have $[X_{-\delta }, X_{-\delta '}]=0$ for the root
vectors $X_{-\delta }, X_{-\delta '}$. (Indeed, otherwise $\delta ,
\delta '\in R^{-}_W$ and $l(w)>2$.)

On the other hand, any map $f_w:\fn_{-\delta}\tto
\fg/(\fg_-\oplus\fl^{(1)})$ which vanishes outside the subspace
$\fn_{-\gamma }\oplus\fn_{-\gamma '}$, where $\{\gamma, \gamma '\}=
R^{-}_W$, and such that the weight of $f_w(\fn_{-\gamma })$ is equal
to $\gamma'$ whereas that of $f_w(\fn_{-\gamma '})$ is $\gamma$
belongs, clearly, to $M_{\rho -w(\rho )}$. This means that $\dim~
M_{\rho -w(\rho )} \leq 2$.

Let $w = s_{\alpha _{1}}s_{\alpha _{2}}$, where $\alpha _{2}\in B_1$
and $(\alpha _{1}, \alpha _{2})\not = 0$. Then $\dim~ M_{\rho
-w(\rho)} = 1$, as the following calculation shows: let
$$
f(X_{-\alpha
_{2}}) = cX'_{\alpha _{1}+\alpha _{2}}, \quad f(X_{-(\alpha _{1}+\alpha
_{2})}) = c'X'_{\alpha _{2}},
$$
where the primed vectors belong to the quotient space. Then, for an
appropriate linear constraint on $c$ and $c'$, we have
\begin{multline*}
X_{-(\alpha _{1}+\alpha _{2})}f(X_{-\alpha _{2}}) -
X_{-\alpha _{2}}f(X'_{ -(\alpha _{1}+\alpha _{2})}) =\\
-cH'_{\alpha _{1}+\alpha _{2}} + c'H'_{ \alpha _{2}} \in\Cee
H'_{\alpha _{1} }\subset(\fg_-\oplus\fl)/(\fg_-\oplus\fl^{(1)}).
\end{multline*}

If $w = s_{\alpha _{1}}s_{\alpha _{2}}$, where $\alpha _{1},
\alpha_{2}\in B_I$, then $\dim~M_{\rho -w(\rho )} = 2$. Since, for any
such $w$, we have
$$
M_{\rho -w(\rho )}\cap Z^2(\fg_-; \fg/(\fg_-\oplus\fl^{(1)})) \neq 0,
$$
we get the desired heading 1)
of the corollary.

The general statement of heading 2) is straightforward; its
particular case a) follows from evident remarks: $H^1(\fg_-; \fg_-^*)
= \fg_-^*\otimes \fg_-^*$ and $H^2(\fg_-) = E^2(\fg_-^*)$.

Proof of heading 2b). Set $\fg_-' = \fg_-/\fz(\fg_-)$.
Clearly, $\fg_-'$ is a trivial $\fg_-'$-module. Therefore, the
sequence
$$
0 \tto \fg_-' \tto \fg_-'\tto \Cee\tto 0
$$
is exact, hence, so is the sequence
$$
0 \tto \Cee\tto \fg_-'\otimes H^1(\fg_-) \tto H^1(\fg_-; \fg_-')
\buildrel f\over\tto H^2(\fg_-) \tto \ldots .
$$
It is not difficult to notice that the image of any cocycle $\fg_-
\tto \fg_-'$ under $f$ belongs to $E^2(\fg_-')$, hence, $f$ is
zero. But then $H^1(\fg_-; \fg_-')$ is isomorphic to the
subspace $(\fg_-'\otimes \fg_-')_0$ of traceless operators
in $\fg_-'\otimes \fg_-'$. It remains to notice that
$H^2(\fg_-)\oplus \Cee = E^2(\fg_-')$. \end{proof}

\ssec{4.5. The number of $\fg _0$-modules} The following Theorem
helps to verify the result. Let $IR$ be the number of irreducible
components in the $\fg _0$-module $H^2(\fg_-; (\fg_-, \fg_0)_*)$.

\begin{Theorem} $IR = \frac{1}{2}c (c + 1) +
\mathop{\sum}\limits c_i$.
\end{Theorem}

\begin{proof} Since $IR = \card ~ W(I)_2$, let us list
the length 2 elements of $W(I)$. Clearly,
$$
W(R_I)_2 = \{w\in W(R_I)\mid l(w) =
2\}\subset W(I)_{ 2}.
$$
If $s_is_j\in W(I)_{2}\setminus W(R_I)_2$, then $s_j\in B_I$ and
$s_i\not \in B_I$. (Indeed, if $s_j\not \in B_I$, then $s_is_j(a_j)
= -s_i(a_j) < 0$ which is false.)

Furthermore, $(\alpha _i, \alpha_j) \not =0$ since otherwise
$s_is_j(\alpha _i) = -\alpha _i < 0$.

Let $\rk ~R = r$. Let $R = R_1\coprod \ldots \coprod R_{s}$, be the
representation as the union of connected components. Let $B =
\{\alpha_1, \ldots , \alpha _r\}$ be a base of $R$. Then
$$
\renewcommand{\arraystretch}{1.4}
\begin{array}{l}
    W(R)_2 = \{s_is_j\; \text{ for
    }\; i\neq j\} =\\
\{s_is_j\; \text{ for }\; i<j\}\coprod
\{s_is_j\; \text{ for }\; i>j\mid
(\alpha _i, \alpha _j) \neq 0\}.
\end{array}
$$
It follows that (here edges are counted ignoring
multiplicities)
\begin{align*}
\card~ W(R)_2 &= \frac{c(c-1)}{2} + \card~\{\text{edges of the Dynkin graph
of } R\}\\
&=\frac{c(c -1)}{2} + c - s.
\end{align*}

If $\alpha \in R\setminus R_I$ and $(\alpha , R_I^{(i)}) \not = 0$
for some $i$, then there exists a unique $\beta\in R_I^{(i)}$ such that
$(\alpha , \beta) \not = 0$. Indeed, if there are two such roots, say,
$\beta _1$ and $\beta _2$, then the Dynkin graph contains a cycle
$$
\alpha\text{ ---}\beta_1\text{ ---}\gamma_1\text{ ---}\gamma_2\text{
---}\ldots\text{ ---} \gamma_s\text{ ---}\beta_2\text{
---}\alpha\text{ ---}
$$
Indeed: $D_I^{(i)}$ is connected and $\alpha \not \in R_I^{(i)}$,
hence, $\gamma _1, \ldots , \gamma _s\in D_I^{(i)}$.

All this demonstrates that
\begin{align*}
\card ~(W_I)_2 &= \card~ (W(R_I))_2 + \sum\limits c_i+1\\
&=\frac{1}{2}c(c -1) + c - s + \sum c_i + s =
\frac{1}{2}c(c -1) + \sum\limits c_i. \qed
\end{align*}
\noqed\end{proof}

\section*{\S 5. The explicit results: the simplest flags}

We consider the standard numbering of vertices of the Dynkin graphs
(same as in \cite{OV} or \cite{Bu}).

Let $k_{i}$ be the coefficient of the $i$th simple root in the
expansion of the maximal root with respect to the simple roots.
For the nilpotent algebra $\fg_{-}=\mathop{\oplus}\limits_{-d \leq
i\leq -1}\fg_{i}$ opposite the maximal parabolic subalgebra
$\fp_{i}$ (with the $i$th selected simple coroot), we have
$d=k_{i}$. The vertices labelled by a $k_{i}=1$ correspond to the
Hermitian symmetric cases already considered in \cite{Go} for the
conformal case and in \cite{LPS} for the Riemannian case. For the
$X=G/P$, we give the lowest weights of the (nonholonomic only for
$k_{i}\neq 1$) N-Weyl tensors with respect to the Levi subgroup
corresponding to the $\Zee$-grading of $\fg$ for which the
selected coroot (vertex) is of degree $-1$, the other coroots
being of degree 0.

The weights of the simple roots are given by the columns of the Cartan
matrix. Set $H^1:=H^1(\fg_-; (\fg_-\oplus \fz)^*)$ and $H^2:=H^2(\fg_-; \fg)$.

\begin{Theorem} Tables $1-4$ give the {\bf lowest} weights of irreducible
$\fl$-components in
$H^2$ in terms of the Cartan matrix (CM): they are
easier to compute, and in terms of the fundamental weights (FW): they
are more conventional.\footnote{ {\bf SuperLie} converts
FW-weights to CM-weights or the other way round in seconds.}

The {\bf highest} weights of irreducible $\fl$-components in $H^1$
are given in terms of the fundamental weights only. To save space,
the column $H^1$ contains all the irreducible $\fl$-modules we
have found but one: for each algebra, for the $i$th node, there
{\bf always} is a component with the highest weight $(0, \ldots ,
0, 2, 0, \ldots , 0)$ with the $2$ on the $i$th place.
\end{Theorem}

We give the weights of the $\widehat{\fg_{0}}=\fl^{(1)}$-modules that
constitute non-conformal part of the N-Riemannian tensor with respect
to the same Cartan subalgebra of $\fg_{0}=\fl$ we used to describe
N-Weyl tensors.

The order of N-Riemannian tensors coincides with the number
occupying the place whose number is equal to the number of the
node. Therefore orders of all structure functions are equal to 2,
except one cases of order 4 for $\fg(2)$ and several cases of
order 3 for $\ff(4)$ and $\fo(2n+1)$.

\begin{Remark} The positive coordinates of lowest weights are not
mistakes: when expressed in terms of the Cartan matrix (CM) they
only occur at the place governed by the center of $\fl$; the
remaining (non-positive) coordinates correspond to the semi-simple
part of $\fl$.
\end{Remark}

\section*{\S 6. The explicit results: case 3) of Yamaguchi's theorem}

\begin{Theorem} In the cases of heading {\em 3)} of Yamaguchi's
theorem, (for completeness we consider the last \lq\lq
selected\rq\rq simple coroot of $\fsl$ as well) the lowest weights
of the $\fg_0$-module $H^2(\fg_-; \fg)$ and highest weights of
$H^1(\fg_-; (\fg_-\oplus \fz)^*)$ and their degrees are listed
below in this section under Tables $1-4$. When the cocycle is
simple-looking it is also given. \end{Theorem}

\noindent \underline{a) $\fg=\fsl(n+1)$}: (the weights are given
with respect to the standard generators of the Cartan subalgebra
of $\fg$, i.e., $e^1_1-e^2_2$, $e^2_2-e^3_3$, \ldots,
$e^n_n-e^{n+1}_{n+1}$)

\noindent \underline{Selected simple coroots: $(1, 2)$}

1) $\deg=0$ of weight $(-3,3,0,-1,0,...,0,-1)$;

2) $\deg=2$ of weight $(4,-1,-1,0,...,0,-1)$;

3) $\deg=3$ of weight $(0,4,-3,0,...,0,-1)$.

\noindent \underline{Selected simple coroots: $(1, i)$, where
$3\leq i\leq n-2$}. To make the answer graphic, the weight is
represented as the sum of a constant part and a part that depends
on $i$ (the $i$th place being $-3$, $-3$ and $2$, respectively;
this variable summand only appears three times):

1) $\deg=0$ of weight $(-1,0,...,0,-1)+(...,0,-1,0,-3,-2,0,...)$;

2) $\deg=0$ of weight $(-1,0,...,0,-1)+(...,0,-2,-3,0,-1,0,...)$;

3) $\deg=1$ of weight $(3,-2,0,...,0,-1)+(...,0,-1,2,-1,0,...)$;

4) $\deg=1$ of weight $(4, -1, -1, 0, 0, 0, -1)$.

\noindent \underline{Selected simple coroots: $(1,n-1)$} (one
weight differs from the above)

1) $\deg=0$ of weight $(-1,0,...,0,-1,0,3,-3)$ (as above);

2) $\deg=1$ of weight $(-1,0,...,0,-3,4,0)$;

3) $\deg=1$ of weight $(3,-2,0,...,0,-1,2,-2)$ (as above);

4) $\deg=1$ of weight $(4,-1,-1,0,...,0,-1)$ (as above).

\noindent \underline{Exceptional cases}:

\noindent \underline{$\fsl(5)$} \underline{Selected simple
coroots: $(1,2)$}

1) $\deg=0$ of weight $(-3, 3, 0, -2)$;

2) $\deg=2$ of weight $(4, -1, -1, -1)$;

3) $\deg=3$ of weight $(0, 4, -3, -1)$.

\noindent \underline{$\fsl(5)$} \underline{Selected simple
coroots: $(1,3)$}

1) $\deg=0$ of weight $(-2, 0, 3, -3)$;

2) $\deg=1$ of weight $(-1, -3, 4, 0)$;

3) $\deg=1$ of weight $(3, -3, 2, -2)$;

4) $\deg=1$ of weight $(4, -1, -1, -1)$.

\noindent \underline{$\fsl(4)$} \underline{Selected simple
coroots: $(1,2)$}

1) $\deg=1$ of weight $(-4,4,0)$;

2) $\deg=2$ of weight $(4,-1,-2)$;

3) $\deg=3$ of weight $(0,4,-4)$.

\noindent \underline{$\fsl(4)$} \underline{Selected simple
coroots: $(1, n)$}:\footnotesize
$$
\renewcommand{\arraystretch}{1.4}
\begin{tabular}{|c|c|c|c|}
\hline
degree&vector: $n>4$&weight: $n>4$&weight: $n=4$\cr
\hline
$1$&$(e^{n+1}_2)^* de^2_1de^3_1$& $(4, -1, -1, 0, ..., 0, -1)$&$(4, -1, -2)$\cr
\hline
$1$&symmetric to the above under&symmetric &symmetric \cr
&$n+1\longleftrightarrow 1$, $n\longleftrightarrow 2$, etc.&to the above&to the
above\cr
\hline
$2$&$(e^{n}_2)^*de^2_1de^{n+1}_{n}$&
$(3, -2, 0, ..., 0, -2, 3)$&$(3, -4, 3)$\cr
\hline
\end{tabular}
$$
\normalsize

\noindent\underline{The exceptional case $\fsl(3)$}:

1) $\deg=4$: $(e^2_3)^*d(e^3_2)d(e^3_1)$ of weight $(-1, 5)$;

2) symmetric to the above.

\noindent\underline{$H^1$}: all vectors except (3) exist for $2\leq i\leq n-1$

(1) $\deg=1$: $((i-1)e^1_1-e^2_2-...-e^i_i)^* de^{i+1}_i$ of weight
 $(...,0,-1,2,-1,.0,...)$;

(2) $\deg=1$: $((i-1)e^1_1-e^2_2-...-e^i_i)^* de^2_1$ of weight
 $(2,-1,0,...)$;

(3) $\deg=2$: $(e^{i+1}_i)^* de^{i+2}_{i-1} - (e^{i+1}_{i-1})^* de^{i+2}_i -
 (e^{i+2}_i)^* de^{i+1}_{i-1} + (e^{i+2}_{i-1})^* de^{i+1}_i$ of
 weight\\
 $(...,0,-1,0,2,0,-1,0,...)$ (for $i=n-1$ of weight
 $(...,0,-1,0,2,0)$), exists for $
 3\leq i\leq n-1$;

(4) $\deg=2$: $(e^{i+1}_i)^* de^{i+1}_i$ of weight
 $(...,0,-2, 4, -2,0,...)$;

(5) $\deg=2$: of weight
 $(1,0,...,0,1,-1,0,...)$:
$$
\renewcommand{\arraystretch}{1.4}
\begin{array}{l}
(n+1)((i-1)e^1_1-e^2_2-...-e^i_i)^* de^{i+1}_1 -\\
(i-1)(n-i+1)\mathop{\sum}\limits_{2\leq j\leq i} ((e^{i+1}_j)^* de^j_1) +
(n+2-i)\mathop{\sum}\limits_{2\leq j\leq i} ((e^j_1)^* de^{i+1}_j);
\end{array}
$$

(6) $\deg=2$: $(e^2_1)^* de^2_1$ of weight
 $(4,-2,0,...)$.

\noindent\underline{b) $\fsp(2n)$}:
\footnotesize
$$
\renewcommand{\arraystretch}{1.4}
\begin{tabular}{|c|c|c|c|}
\hline
degree&weight: $n>3$&weight: $n=3$&$H^1$\cr
\hline
$-1$&$(-2, 0, \ldots , 0, -1, -2, 3)$&$(-3, -2, 3)$&$-$\cr
\hline
$1$&$(4, -3, 0, \ldots , 0, -2, 2)$&$(4, -5, 2)$&$(0, \ldots, 0, -2, 2)$\cr
$1$&$(5, -2, -1, 0, \ldots , 0)$&$(5, -2, -1)$&$(2, -1, 0, \ldots, 0)$\cr
\hline
\end{tabular}
$$\normalsize
There are also irreducible components in $H^1$ of degree 2 and their weights are
$$
(0, \ldots, 0, -2, 0, 2), \quad (0, \ldots, 0, -4, 4), \quad (1,
0, \ldots , 0, -1, 1), \quad (4, -2, 0, \ldots , 0).
$$

\noindent\underline{The exceptional case $\fsp(4)$}:

\underline{$H^2$}: $\deg=3$ of weight $(6, -3)$ and $\deg=4$ of weight $(-4,
5)$.

\underline{$H^1$:} $\deg=1$ of weights $(-2, 2)$ and $(2, -1)$
$\deg=2$ of weights $(-4, 4)$, $(0, 1)$ and $(4, -2)$

\section*{\S 7. The explicit results: the $G(2)$-structures}

\ssec{7.1. The $G(2)$-structure} Let $M$ be a manifold with the
$G=G(2)$-structure, i.e., the $G(2)$-module $T_{m}M$ is isomorphic
to the 1st fundamental module, cf. \cite{B}. We wondered to what
extent the \lq\lq positive 3-form" (an invariant of the
$G(2)$-structure similar to the metric in the Riemannian case,
i.e., for the $\fo(n)$-structure), see \cite{B}, can be reduced to
a canonical form. Below is a description of the space of
obstructions to canonicity. Since $\fg(2)\subset\fo(7)$ and
$\fo(n)_{1}=0$ for $n>2$, it is easy to see that $(\fg_{-1},
\fg_{0})_*=\fg_{-1}\oplus \fg_{0}$.

\begin{Statement} As $\fg_{0}$-module, $H^2(\fg_{-1}; (\fg_{-1},
\fg_{0})_*)$ is the direct sum of the irreducible $G(2)$-modules
whose highest weights and orders of the corresponding structure
function are given by the following table:
$$
\renewcommand{\arraystretch}{1.4}
\begin{tabular}{|c|c|}
\hline
weight&order\cr
\hline
$(0, 0), \;(0, 1), \;(1, 0), \;(2, 0)$&$1$\cr
\hline
$(0, 2)$&$2$\cr
\hline
\end{tabular}
$$
\end{Statement}
The requirement of vanishing of order 1 structure functions (for the
corresponding equations, see \cite{B}) is an analog of Wess-Zumino
constraints in supergravity \cite{WB}.

From the list of $\Zee$-gradings of simple Lie algebras we know
that there is no analog of classical domain with the
$G=CG(2)$-structure, i.e., conformal, even nonholonomic one.
Contrariwise, one and only one of the exceptional Lie
superalgebras has such a grading. Let us compute the corresponding
space of structure functions.

\ssec{7.2. The $CG(2)$-structure on $\Cee P^{1, 7}$ with a
nonholonomic distribution} We consider $\Cee P^{1, 7}$ as the
quotient of the simple Lie supergroup $AG(2)$ modulo the parabolic
subalgebra corresponding to the grading $(1, 0, 0)$ for the Cartan
matrix, cf. \cite{GL2}, where the Lie superalgebra
$\fag(2)=\text{Lie}(AG(2))$ is presented:
$$
\begin{pmatrix}
0&1&0\\
-1&2&-3\\
0&-1&2
\end{pmatrix}.
$$
Here $\fg_{0}=\fg(2)\oplus\fz$ for a 1-dimensional $\fz$; let
$\widehat{\fg_{0}}=\fg(2)$. It is now not as easy as in sec. 7.1
to see that $(\fg_{-}, \fg_{0})_*=\fag(2)$ and $(\fg_{-},
\widehat{\fg_{0}})_*=\fg_{-}\oplus \widehat{\fg_{0}}$, but still
true. No super version of the BWB theorem exists (cf. \cite{Pe},
\cite{PS}) to help us, so, to obtain the following statement, we
used {\bf SuperLie}.

\begin{Statement} As $\fg_{0}=\fsl(2)\oplus \fg(2)$-module,
$H^2(\fg_{-1}; \fag(2))$ --- the space of N-conformal structure
functions --- is an irreducible module with highest weight $(5, 0,
1)$.

As $\widehat{\fg_{0}}=\fz\oplus \fg(2)$-module, $H^2(\fg_{-1};
(\fg_{-}, \widehat{\fg_{0}})_*)$ --- the space of N-Riemannian
structure functions --- is the direct sum of irreducible modules
whose highest weights are as follows:
$$
(5, 0, 1), \quad (6, 1, 0), \quad (7, 0, 1).
$$
\end{Statement}

\begin{Remark} In our attempt to understand the meaning of tensors $V(\fg_i)$ in
\cite{B}, where {\it exterior} powers of $\fg_{-1}$ appear, we
conjectured that these $V(\fg_i)$ might be related with the
$G(2)$-structure on a {\it purely odd}\ superspace. However,
having computed the corresponding structure functions

\footnotesize
$$
\renewcommand{\arraystretch}{1.4}
\begin{tabular}{|c|c|}
\hline
weight&order\cr
\hline
$(0, 1), \;(1, 0), \;(3, 0)$&$1$\cr
\hline
$(0, 1)$&$2$\cr
\hline
\end{tabular}
$$\normalsize
we see that they do not coincide with the $V(\fg_i)$, so the
meaning of these $V(\fg_i)$'s remains a mystery to us.
\end{Remark}

\vfill
\newpage

\tiny
\vbox{
$$
\renewcommand{\arraystretch}{1.4}
\begin{array}{l}
\text{Table 1}\\
\\
\renewcommand{\arraystretch}{1.4}
\begin{tabular}{|c|c|c|c|c|c|}
\hline $\fg$&Node&$k_{i}$&CM&FW&$H^1$\cr \hline
$\fg(2)$&$1$&$3$&$(8, -4)$&$(4, 0)$&$(4, 2)$\cr \cline{2-6}
&$2$&$2$&$(-7, 4)$&$(-2, 1)$&$(2, 2)$\cr \hline
$\ff(4)$&$1$&$2$&$(3, 0, -1, -1)$&$(0, -3, -3, -2)$&$(2, 3, 2, 1)$
\cr \cline{2-6} &$2$&$4$&$(0, 3, -2, -1)$&$(-1, -2, -3, -2)$&$(0,
3, 2, 1)$ \cr &&&$(-3, 4, -1, -2)$&$(-2, -1, 2, -2)$&$(1, 2, 1,
0)$\cr \cline{2-6} &$3$&$3$&$(0, -6, 4, 0)$&$(-2, -4, 0, 0)$&$(0,
2, 2, 0)$ \cr &&&$(-1, -2, 3, -3)$&$(-2, -3, -1, -2)$&$(0, 1, 2,
0)$\cr \cline{2-6} &$4$&$2$&$(0, -2, -1, 4)$&$(-2, -4, -2,
1)$&$(0, 2, 2, 2)$\cr \hline $\fe(6)$&$1$&$1$&$(3, 0, -1, 0, 0,
-1)$&$(1, -1, -3, -2, -1, -2)$&$(2, 2, 2, 1, 0, 1)$\cr \cline{2-6}
&$2$&$2$&$(0, 3, -2, 0, 0, -1)$&$(0, 0, -3, -2, -1, -2)$& $(1, 2,
1, 0, 0, 0)$\cr &&&$(-2, 3, 0, -1, 0, -2)$&$(-1, 0, -2, -2, -1,
-2)$&$(0, 2, 2, 1, 0, 1)$\cr \cline{2-6} &$3$&$3$&$(0, -3, 4, -3,
0, 0)$&$(-1, -2, 0, -2, -1, 0)$& $(0, 1, 2, 1, 0, 0)$\cr &&&$(0,
-2, 3, 0, -1, -3)$&$(-1, -2, -1, -1, -1, -2)$&$(0, 0, 2, 1, 0,
1)$\cr &&&$(-1, 0, 3, -2, 0, -3)$&$(-1, -1, -1, -2, -1, -2)$&$(0,
1, 2, 0, 0, 1)$\cr \cline{2-6} &$4$&$2$&$(0, 0, -2, 3, 0,
-1)$&$(-1, -2, -3, 0, 0, -2)$& $(0, 0, 1, 2, 1, 0)$\cr &&&$(0, -1,
0, 3, -2, -2)$&$(-1, -2, -2, 0, -1, -2)$& $(0, 1, 2, 2, 0, 1)$\cr
\cline{2-6} &$5$&$1$&$(0, 0, -1, 0, 3, -1)$&$(-1, -2, -3, -1, 1,
-2)$& $(0, 1, 2, 2, 2, 1)$\cr \cline{2-6} &$6$&$2$&$(0, -1, -1,
-1, 0, 4)$&$(-1, -2, -2, -2, -1, 1)$& $(0, 1, 2, 1, 0, 2)$\cr
\hline $\fe(7)$&$1$&$1$&$(3, 0, -1, 0, 0, -1, 0)$& $(1, -1, -3,
-4, -3, -2, -2)$&$(2, 2, 2, 2, 1, 0, 1)$\cr \cline{2-6}
&$2$&$2$&$(0, 3, -2, 0, 0, -1, 0)$&$(0, 0, -3, -4, -3, -2, -2)$&
$(1, 2, 1, 0, 0, 0, 0)$\cr &&&$(-2, 3, 0, -1, 0, -1, 0)$&$(-1, 0,
-2, -4, -3, -2, -2)$&$(0, 2, 2, 2, 1, 0, 1)$\cr \cline{2-6}
&$3$&$3$&$(0, -2, 3, 0, -1, -1, -1)$&$(-1, -2, -1, -3, -3, -2,
-2)$& $(0, 1, 2, 1, 0, 0, 0)$\cr &&&$(-1, 0, 3, -2, 0, -1,
0)$&$(-1, -1, -1, -4, -3, -2, -2)$&$(0, 0, 2, 2, 1, 0, 1)$\cr
\cline{2-6} &$4$&$4$&$(0, 0, -2, 3, -2, -1, 0)$&$(-1, -2, -3, -2,
-3, -2, -1)$& $(0, 0, 1, 2, 1, 0, 0)$ \cr &&&$(0, 0, -2, 3, 0, -2,
-2)$&$(-1, -2, -3, -2, -2, -2, -2)$&$(0, 0, 0, 2, 1, 0, 1)$\cr
&&&$(0, -1, 0, 3, -2, -1, -2)$&$(-1, -2, -2, -2, -3, -2, -2)$&$(0,
0,
 1,
2, 0, 0, 1)$\cr
\cline{2-6}
&$5$&$3$&$(0, 0, 0, -3, 4, 0, 0)$&
$(-1, -2, -3, -4, 0, 0, -2)$&
$(0, 0, 0, 1, 2, 1, 0)$\cr
&&&$(0, 0, -1, 0, 3, -3, -1)$&
$(-1, -2, -3, -3, -1, -2, -2)$&$(0, 0, 1, 2, 2, 0, 1)$\cr
\cline{2-6}
&$6$&$2$&$(0, 0, 0, -1, -1, 4, 0)$&$(-1, -2, -3, -4, -2, 1, -2)$&
$(0, 0, 1, 2, 2, 2, 1)$\cr
\cline{2-6}
&$7$&$2$&$(0, 0, -1, 0, -1, -1, 3)$&$(-1, -2, -3, -3, -3, -2, 0)$&
$(0, 0, 1, 2, 1, 0, 2)$\cr
\hline
$\fe(8)$&$1$&$2$& $(4, -1, -1, 0, 0, 0, 0, 0)$&$(-1, -2, -4, -5,
-6, -4, -2, -3)$& $(2, 2, 2, 2, 2, 1, 0, 1)$\cr \cline{2-6}
&$2$&$3$&$ (0, 4, -3, 0, 0, 0, 0, 0)$& $(0, 0, -4, -5, -6, -4, -2,
-3)$& $(1, 2, 1, 0, 0, 0, 0, 0)$ \cr &&&$(-3, 3, 0, -1, 0, 0, 0,
0)$& $(-2, -1, -3, -5, -6, -4, -2, -3)$&$(0, 2, 2, 2, 2, 1, 0,
1)$\cr \cline{2-6} &$3$&$4$&$(-1, -2, 3, 0, -1, 0, 0, 0)$& $(-2,
-3, -2, -4, -6, -4, -2, -3)$& $(1, 2, 1, 0, 0, 0, 0, 0)$ \cr
&&&$(-2, 0, 3, -2, 0, 0, 0, 0)$& $(-2, -2, -2, -5, -6, -4, -2,
-3)$&$(0, 2, 2, 2, 2, 1, 0, 1)$\cr \cline{2-6} &$4$&$5$&$(-1, 0,
-2, 3, 0, -1, 0, -1)$& $(-2, -3, -4, -3, -5, -4, -2, -3)$& $(0, 0,
1, 2, 1, 0, 0, 0)$ \cr &&&$(-1, -1, 0, 3, -2, 0, 0, 0)$& $(-2, -3,
-3, -3, -6, -4, -2, -3)$&$(0, 0, 0, 2, 2, 1, 0, 1)$\cr \cline{2-6}
&$5$&$6$& $(-1, 0, 0, -2, 3, -2, 0, 0)$&$(-2, -3, -4, -5, -4, -4,
-2, -2)$& $(0, 0, 0, 1, 2, 1, 0, 0)$ \cr &&&$(-1, 0, 0, -2, 3, 0,
-1, -2)$&$(-2, -3, -4, -5, -4, -3, -2, -3)$&$(0, 0, 0, 0, 2, 1, 0,
1)$\cr &&&$(-1, 0, -1, 0, 3, -2, 0, -2)$&$(-2, -3, -4, -4, -4, -4,
-2, -2)$&$(0, 0, 0, 1, 2, 0, 0, 1)$\cr \cline{2-6} &$6$&$4$&$(-1,
0, 0, 0, -2, 3, 0, 0)$& $(-2, -3, -4, -5, -6, -3, 0, -3)$& $(0, 0,
0, 0, 1, 2, 1, 0)$ \cr &&&$(-1, 0, 0, -1, 0, 3, -2, -1)$&$(-2, -3,
-4, -5, -5, -4, -2, -1)$&$(0, 0, 0, 1, 2, 2, 0, 1)$\cr \cline{2-6}
&$7$&$2$&$(-1, 0, 0, 0, 0, -1, 2, 0)$&$(-1, 0, 0, 0, -1, 0, 3,
0)$& $(0, 0, 0, 1, 2, 2, 2, 1)$\cr \cline{2-6} &$8$&$3$&$(-1, 0,
0, 0, -1, 0, 0, 2)$& $(-1, 0, 0, -1, 0, -1, 0, 3)$& $(0, 0, 0, 1,
2, 1, 0, 2)$\cr \hline
\end{tabular}
\end{array}
$$
}

\vbox{
\tiny
$$
\renewcommand{\arraystretch}{1.4}
\begin{array}{l}
\text{Table 2 } \fo(2n):\\
\\
\renewcommand{\arraystretch}{1.4}
\begin{tabular}{|c|c|c|c|c|}
\hline
Node&$k(i)$&$rk$&weight&$H^1$\cr
\hline
$1$&$1$&$4$&$(2, 0, -1, -1)$&$(2211)$\cr
\cline{3-5}
&&$4+l$&$(2, 0, \underbrace{-2, \ldots, -2}_{l}, -1, -1)$&
$(\underbrace{2 \ldots 2}_{l-2}, 1 1)$\cr
\hline
$2$&$2$&$5+l$&$(0, 1, -2, \underbrace{-2, \ldots, -2}_{l}, -1,
-1)$&$(02211)$\cr
&&&$(-1, 1, -1, \underbrace{-2, \ldots, -2}_{l}, -1, -1)$&$(12100)$\cr
\hline
$3$&$2$&$6+l$&$(-1, -2, 0, 1, \underbrace{-2, \ldots, -2}_{l}, -1, -1)$&
$(0 \underbrace{2 \ldots 2}_{l-3} 1 1)$\cr
&&&$-1, 0, 1, -2, \underbrace{-2, \ldots, -2}_{l}, -1, -1)$&
$(1 2 1 \underbrace{0 \ldots 0}_{l-3})$\cr
\hline
$3+k$&$2$&$7$&$(-1, \underbrace{-2, \ldots, -2}_{k} -2, 0, -1, -1,
-1)$&$(002211)$\cr
$k<rk-6$&&&$(-1, \underbrace{-2, \ldots, -2}_{k}, -1, 0, -2, -1,
-1)$&$(012100)$\cr
\cline{3-5}
&&$>7$&$(-1, \underbrace{-2, \ldots, -2}_{k} -2, 0, -1,
\underbrace{-2, \ldots, -2}_{l}, -1, -1)$&$(0 \ldots 022 1 1)$\cr
&&&$(-1, \underbrace{-2, \ldots, -2}_{k}, -1, 0, -2,
\underbrace{-2, \ldots, -2}_{l}, -1, -1)$&$(0 \ldots 0121 00)$\cr
\hline
$\dots$ & $\dots$ & $\dots$ &$\dots$&$\dots$\cr
\hline
fork&2&4&$(0, 1, -1, -1)$&$(0211)$ \cr
&&&$(-1, 1, -1, 0)$&$(1201)$\cr
&&&$(-1, 1, 0, -1)$&$(1210)$\cr
\cline{3-5}
&&5&$(-1, -2, 0, -1, 0)$&$(00211)$\cr
&&&$(-1, -2, 0, 0, -1)$&$(01201)$\cr
&&&$(-1, 0, 1, -1, -1)$&$(01210)$\cr
\cline{3-5}
&&$5+k$&$(-1, \underbrace{-2, \ldots, -2}_{k}, -2, 0, -1, 0)$&$(0\ldots
00211)$\cr
&&$k<rk-5$&$(-1, \underbrace{-2, \ldots, -2}_{k}, -2, 0, 0, -1)$&$(0\ldots
01201)$\cr
&&&$(-1, \underbrace{-2, \ldots, -2}_{k}, -1, 0, -1, -1)$&$(0\ldots
01210)$\cr
\hline
\end{tabular}\\
\end{array}
$$
}

\vbox{
\tiny
$$
\renewcommand{\arraystretch}{1.4}
\begin{array}{l}
\text{Table 3 } \fo(2n+1):\\
\\
\renewcommand{\arraystretch}{1.4}
\begin{tabular}{|c|c|c|c|c|}
\hline
Node&$k(i)$&$rk$&weight&$H^1$\cr
\hline
1&1&2&$(3, 1)$&$(22)$\cr
\cline{3-5}
&&3&$(2, 0, -2)$&$(222)$\cr
\cline{3-5}
&&$3+k$&$(2, 0, \underbrace{-2, \ldots, -2}_{k})$&$(2\ldots 2)$\cr
\hline
2&2&3&$(0, 1, -2)$&$(022)$\cr
&&&$(-1, 1, -1)$&$(121)$\cr
\cline{3-5}
&&$3+k$&$(0, 1, -2, \underbrace{-2, \ldots, -2}_{k})$&$(02\ldots 2)$\cr
&&&$(-1, 1, -1, \underbrace{-2, \ldots, -2}_{k})$&$(1210\ldots 0)$\cr
\hline
3&2&$4+k$&$(-1, -2, 0, 1, \underbrace{-2, \ldots, -2}_{k})$&$(002\ldots
2)$\cr
&&&$(-1, 0, 1, -2 \underbrace{-2, \ldots, -2}_{k})$&$(01210\ldots 0)$\cr
\hline
$\dots$ & $\dots$ & $\dots$ &$\dots$&$\dots$\cr
\hline
penul-&2&$5+k$&$(-1, \underbrace{-2, \ldots, -2}_{k} -2, 0, -1,
-2)$&$(0\ldots 022)$\cr
timate&&&$(-1, \underbrace{-2, \ldots, -2}_{k}, -1, 0, -2, -2)$&$(0\ldots
0121)$\cr
\hline
last&1&2&$(0, 3)$&$-$\cr
\cline{3-5}
&&3&$(-1, 0, 3)$&$(123)$\cr
\cline{3-5}
&&$3+k$&$(-1, \underbrace{-2, \ldots, -2}_{k}, -1, 1)$&$(0\ldots 0123)$\cr
\hline
\end{tabular}\\
\end{array}
$$
}

\tiny 
$$
\renewcommand{\arraystretch}{1.4}
\begin{array}{l}
\text{Table 4 } \fsp(2n):\\
\\
\renewcommand{\arraystretch}{1.4}
\begin{tabular}{|c|c|c|c|c|}
\hline
Node&$k(i)$&$rk$&weight&$H^1$\cr
\hline
1&2&2&$(3, 0)$&$-$\cr
\cline{3-5}
&&3&$(2, 1, -1)$&$-$\cr
\cline{3-5}
&&$3+k$&$(2, -1, \underbrace{-2, \ldots, -2}_{k}, -1)$&$-$\cr
\hline
2&2&3&$(1, 2, -1)$&$(121)$\cr
&&&$(-2, 1, 0)$&\cr
\cline{3-5}
&&$4+k$&$(1, 2, -2, \underbrace{-2, \ldots, -2}_{k}, -1)$&$(1210\ldots
0)$\cr
&&&$(-2, 0, -1, \underbrace{-2, \ldots, -2}_{k}, -1)$&$(12\ldots 21)$\cr
\hline
$\dots$ & $\dots$ & $\dots$ &$\dots$&$\dots$\cr
\hline
pen-penul-&2&$4+k$&$(\underbrace{-2, \ldots, -2}_{k}, -2, 0, 1,
-1)$&$(0\ldots 0121 0)$\cr
timate&&&$(\underbrace{-2, \ldots, -2}_{k}, 0, 1, -2, -1)$&$(0\ldots
01211)$\cr
\hline
penul-&2&$3+k$&$(\underbrace{-2, \ldots, -2}_{k} -2, 1, 0)$&$(0\ldots
0121)$\cr
timate&&&$(\underbrace{-2, \ldots, -2}_{k}, -1, 0, -1)$&\cr
\hline
last&1&2&$(1, 3)$&$(22)$\cr
\cline{3-5}
&&$3+k$&$(\underbrace{-2, \ldots, -2}_{k}, -2, -1, 1)$&$(0\ldots 022)$\cr
\hline
\end{tabular}\\
\end{array}
$$

\end{document}